\documentclass{amsart}
\usepackage{amsthm}
\usepackage{amssymb}
\usepackage{upref}
\usepackage{amscd}

\numberwithin{equation}{section}
\newtheorem{main}{Theorem}

\newtheorem{theorem}[equation]{Theorem}
\newtheorem{lemma}[equation]{Lemma}
\newtheorem{proposition}[equation]{Proposition}
\newtheorem{corollary}[equation]{Corollary}

\theoremstyle{definition}
\newtheorem{definition}[equation]{Definition}
\newtheorem{remark}[equation]{Remark}

\DeclareMathOperator{\Hom}{Hom}
\DeclareMathOperator{\Ann}{Ann}
\DeclareMathOperator{\sho}{ho}
\DeclareMathOperator{\Ext}{Ext}
\DeclareMathOperator{\Tor}{Tor}
\DeclareMathOperator{\Spec}{Spec}

\DeclareMathOperator{\colim}{colim}

\renewcommand{\smash}{\wedge}

\newcommand{\Q}{\mathbb{Q}}
\newcommand{\ideal}[1]{\mathfrak{#1}}

\newcommand{\Ch}[1]{\text{Ch} (#1)}
\newcommand{\stable}[1]{\text{Stable} (#1)}
\newcommand{\hopfalg}{(A,\Gamma )}

\newcommand{\cat}[1]{\mathcal{#1}}

\newcommand{\Aff}{\mathbf{Aff}}

\newcommand{\Z}{\mathbb{Z}}
\newcommand{\Fp}{\mathbb{F}_{p}}

\newcommand{\comod}{\text{-comod}}
\newcommand{\Mod}{\text{-mod}}

\newcommand{\Bousfield}[1]{<\! #1 \! >}

\newcommand{\mathcolon}{\colon\,}

\newcommand{\uc}{\textup{:}}
\newcommand{\usc}{\textup{;}}

\begin{document}
 
\title{Chromatic phenomena in the algebra of $BP_{*}BP$-comodules} 

\date{\today}

\author{Mark Hovey}
\address{Department of Mathematics \\ Wesleyan University
\\ Middletown, CT 06459}
\email{hovey@member.ams.org}


\begin{abstract}
We describe the author's research with Neil Strickland on the global
algebra and global homological algebra of the category of
$BP_{*}BP$-comodules.  We show,
following~\cite{hovey-strickland-comodules}, that the category of
$E(n)_{*}E(n)$-comodules is a localization, in the abelian sense, of the
category of $BP_{*}BP$-comodules.  This gives analogues of the usual
structure theorems, such as the Landweber filtration theorem, for
$E(n)_{*}E(n)$-comodules.  We recall the work
of~\cite{hovey-comodule-homotopy}, where an improved version
$\stable{\Gamma}$ of the derived category of comodules over a
well-behaved Hopf algebroid $(A, \Gamma)$ is constructed.  The main new
result of the paper is that $\stable{E(n)_{*}E(n)}$ is a Bousfield
localization of $\stable{BP_{*}BP}$, in analogy to the abelian
case.  
\end{abstract}

\maketitle

\section*{Introduction}

The object of this paper is to describe some of the author's recent
work, much of it joint with Neil Strickland, on comodules over
$BP_{*}BP$ and related Hopf algebroids.  The basic idea of this work is
to realize the chromatic approach to stable homotopy theory in the
algebraic world of comodules.  This means, in particular, constructing
and understanding the localization functor $L_{n}$, or the associated
finite localization functor $L_{n}^{f}$, in the abelian category
$BP_{*}BP\comod$ of graded $BP_{*}BP$-comodules.  It also means
constructing $L_{n}$ and $L_{n}^{f}$ in some kind of associated derived
category $\stable{BP_{*}BP}$ of chain complexes of $BP_{*}BP$-comodules.
Ultimately, we would like to understand $L_{K(n)}$ in the algebraic
setting as well, and this should involve the Morava stabilizer groups.
In the present paper, we confine ourselves to $L_{n}$ and $L_{n}^{f}$.

Topologically, $L_{n}^{f}$ is localization away from a finite spectrum
of type $n+1$, and $L_{n}$ is localization with respect to the homology
theory $E(n)$.  These functors are probably different in the ordinary
stable homotopy category (because Ravenel's telescope
conjecture~\cite{rav-loc} is widely expected to be false), but they turn
out to agree on the abelian category of $BP_{*}BP$-comodules.  We then
have the following theorem.

\begin{main}\label{main-A}
Let $L_{n}$ denote localization away from $BP_{*}/I_{n+1}$ in the
category of $BP_{*}BP$-comodules.  Then there is an equivalence of
categories between $E(n)_{*}E(n)$-comodules and $L_{n}$-local
$BP_{*}BP$-comodules.
\end{main}

Note that this is localization in an abelian sense.  The localization
$L_{n}$ turns out to be the localization functor that inverts all
maps whose kernel and cokernel are $v_{n}$-torsion.  

This theorem is really a special case of a more general theorem proved
in~\cite{hovey-strickland-comodules}.  In general, if $\hopfalg$ is a
flat Hopf algebroid and $B$ is a Landweber exact $A$-algebra, there is
an induced flat Hopf algebroid $(B, \Gamma_{B})$.  We prove
in~\cite{hovey-strickland-comodules} that, in this situation, the
category of $\Gamma_{B}$-comodules is always equivalent to some
localization of the category of $\Gamma$-comodules.  In the case of
$BP_{*}BP$, we give a partial classification of such localizations.  In
particular, $E(n)_{*}E(n)$ can be replaced by $E_{*}E$ for any Landweber
exact commutative ring spectrum $E$ with $E_{*}/I_{n}\neq 0$ and
$E_{*}/I_{n+1}=0$.  This tells us that all such theories $E$ have
equivalent categories of $E_{*}E$-comodules, even though the categories
of $E_{*}$-modules can be drastically different.  It also leads to
structural results about $E_{*}E$-comodules analogous to those of
Landweber~\cite{land-exact} for $BP_{*}BP$-comodules.

To extend this result to the derived category setting, we first must
decide what we mean by the derived category.  The ordinary derived
category, obtained by inverting homology isomorphisms, is usually badly
behaved.  For example, the analogue of $S^{0}$ in the derived category
of $BP_{*}BP$-comodules is $BP_{*}$ thought of as a complex concentrated
in degree $0$, but this is not a small object (see the introduction to
Section~\ref{sec-derived}).

The following is a corollary of the main result
of~\cite{hovey-comodule-homotopy}.  

\begin{main}\label{main-B}
Suppose $E$ is a commutative ring spectrum that is Landweber exact over
$BP$.  Then there is a bigraded monogenic stable homotopy category
$\stable{E_{*}E}$ such that 
\[
\pi_{**}S^{0} \cong \Ext^{**}_{E_{*}E}(E_{*}, E_{*}).
\]
\end{main}

In particular, the sphere, which is $E_{*}$ concentrated in degree $0$,
is a small object of $\stable{E_{*}E}$.  This stable homotopy category
is the homotopy category of a suitable model structure on chain
complexes of $E_{*}E$-comodules.  The weak equivalences are homotopy
isomorphisms, where homotopy is suitably defined.  Every cofibrant
object is dimensionwise projective over $E_{*}$ and every complex of
relatively injective comodules is fibrant.  Just as in the ordinary
stable homotopy category, there exist interesting nontrivial complexes
with no homology.

On $\stable{BP_{*}BP}$, we define $L_{n}^{f}$ to be finite localization
away from $BP_{*}/I_{n+1}$, which turns out to be a small object in
$\stable{BP_{*}BP}$.  We define $L_{n}$ to be Bousfield localization
with respect to the homology theory corresponding to $E(n)_{*}$.  In
$\stable{BP_{*}BP}$, the homology functor $H$ is represented by
$BP_{*}BP$, so is somewhat analogous to the $BP$-homology of a spectrum.
Then the homology theory corresponding to $E(n)_{*}$ is in fact $HE(n)$,
ordinary homology with coefficients in $E(n)_{*}$.  Because of this,
some of the things one would expect to be true about $L_{n}$ are false.
For example, the two localizations $L_{n}^{f}$ and $L_{n}$ on
$\stable{BP_{*}BP}$ are definitely different in general, so the most
naive version of the telescope conjecture is false in
$\stable{BP_{*}BP}$.  Also, $L_{n}$ is not a smashing localization,
though $L_{n}^{f}$ is.  Thus $L_{n}$ is less important in
$\stable{BP_{*}BP}$ then $L_{n}^{f}$.  

These statements are also true in the category $\stable{E(n)_{*}E(n)}$,
where $L_{n}^{f}$ is the identity functor (as it is localization away
from $E(n)_{*}/I_{n+1}=0$), and $L_{n}$ is localization with respect to
ordinary homology (which already has coefficients in $E(n)_{*}$).  Thus
$L_{n}\stable{E(n)_{*}E(n)}$ is the classical unbounded derived category
of $E(n)_{*}E(n)$-comodules.

The main new result of this paper is then the following theorem.  

\begin{main}\label{main-C}
There is an equivalence of stable homotopy categories between the
localization $L_{n}^{f}\stable{BP_{*}BP}$ and $\stable{E(n)_{*}E(n)}$.
\end{main}

As above in Theorem~\ref{main-A}, we can replace $E(n)$ in
Theorem~\ref{main-C} by any Landweber exact commutative ring spectrum
$E$ with $E_{*}/I_{n}\neq 0$ and $E_{*}/I_{n+1}=0$.
Theorem~\ref{main-C} also implies a similar equivalence between
$L_{n}\stable{BP_{*}BP}$ and $L_{n}\stable{E(n)_{*}E(n)}$.  

From a computational point of view, Theorem~\ref{main-C} gives rise to
the following change of rings theorem.  

\begin{main}\label{main-D}
Suppose that $M$ is a finitely presented $BP_{*}BP$-comodule and $N$ is
a $BP_{*}BP$-comodule such that $N=L_{n}N$ and the right derived functors
$L_{n}^{i}N$ are $0$ for $i>0$.  Then there is a change of rings
isomorphism 
\[
\Ext_{BP_{*}BP}^{**}(M,N) \cong
\Ext_{E(n)_{*}E(n)}(E(n)_{*}\otimes_{BP_{*}}M,
E(n)_{*}\otimes_{BP_{*}}N).  
\]
\end{main}

This change of rings theorem includes the Miller-Ravenel change of rings
theorem~\cite{miller-ravenel} and the change of rings theorem of the
author and Sadofsky~\cite{hovey-sadofsky-picard} as special cases.
Also, $E(n)$ can be replaced by any Landweber exact commutative ring
spectrum $E$ with $E_{*}/I_{n}\neq 0$ and $E_{*}/I_{n+1}=0$.  

Because the basic structure of the abelian category of comodules over a
flat Hopf algebroid is not as well known as it should be, we first
summarize this in Section~\ref{sec-comod}.  The results in this section
were mostly proved in~\cite{hovey-comodule-homotopy}.  We then describe
our proof of Theorem~\ref{main-A} and related results about
$E(n)_{*}E(n)$-comodules in Section~\ref{sec-Landweber}.  Further
details can be found in~\cite{hovey-strickland-comodules}.  We introduce
the stable homotopy category of comodules in Section~\ref{sec-derived},
describing the proof of Theorem~\ref{main-B} and looking at some
particular features of $\stable{BP_{*}BP}$ and $\stable{E(n)_{*}E(n)}$.
Some of the material in this section can be found
in~\cite{hovey-comodule-homotopy}, but some of it is new.  We discuss
the relation between $\stable{BP_{*}BP}$ and $\stable{E(n)_{*}E(n)}$ in
Section~\ref{sec-comparison}, where we prove Theorems~\ref{main-C}
and~\ref{main-D}.  All of the results in this section are new.

It is a pleasure to thank my coauthor Neil Strickland.  Many of the
theorems in this paper are joint work with him; joint work that started
in Barcelona the week before the 2002 Barcelona conference on algebraic
topology.  I would like to thank the Universitat de Barcelona, the
Universitat Aut\`{o}noma de Barcelona, the Centre de Recerca Matematica,
and the Isaac Newton Institute for Mathematical Sciences for their
support during our collaboration.  I would also like to thank John
Greenlees for many discussions about the material in this paper.

\section{Comodules}\label{sec-comod}

The object of this section is to give an overview of the structural
properties of the category $\Gamma \comod$.  We begin by recalling the
structure maps of a Hopf algebroid $\hopfalg$.

A Hopf algebroid $\hopfalg$ has the following structure maps, which are
all maps of commutative rings. 
\begin{itemize}
\item The \emph{counit} $\epsilon \mathcolon \Gamma \xrightarrow{} A$,
corepresenting the identity map of an object. 
\item The \emph{left unit} $\eta_{L}\mathcolon A\xrightarrow{} \Gamma$,
corepresenting the source of a morphism.  
\item The \emph{right unit} $\eta_{R}\mathcolon A\xrightarrow{} \Gamma$,
corepresenting the target of a morphism.  
\item The \emph{diagonal} $\Delta \mathcolon \Gamma \xrightarrow{}
\Gamma \otimes_{A} \Gamma $, corepresenting the composite of two
composable morphisms.  Note that this is a tensor product of
$A$-bimodules, with the left $A$-module structure given by $\eta_{L}$
and the right $A$-module structure given by $\eta_{R}$. 
\item The \emph{conjugation} $\chi \mathcolon \Gamma
\xrightarrow{}\Gamma$, corepresenting the inverse of a morphism.  
\end{itemize}

There are many relations between these structure maps, but they are all
easily obtained from corresponding facts about groupoids.  For example,
the source and target of the identity morphism at $x$ are both $x$, so
$\epsilon \eta_{L}=\epsilon \eta_{R}=1$.  

Note that the conjugation is best thought of as a map $\chi \mathcolon
\Gamma \xrightarrow{} \widetilde{\Gamma}$, where $\widetilde{\Gamma}$
denotes $\Gamma$ with the opposite $A$-bimodule structure, so that $A$
acts on the left on $\widetilde{\Gamma}$ by $\eta_{R}$.  Similarly, the
multiplication map is best thought of as $\mu \mathcolon \Gamma
\otimes_{A}\widetilde{\Gamma} \xrightarrow{}\Gamma$, although we could
also think of it as having domain $\widetilde{\Gamma}\otimes_{A}\Gamma$.
With this convention, the fact that composition of a map with its
inverse gives the identity gives the following commutative diagram,
\[
\begin{CD}
\Gamma  @>\Delta >> \Gamma \otimes_{A} \Gamma @>1\otimes \chi >> \Gamma
\otimes_{A} \widetilde{\Gamma} \\
@| @. @VV\mu V \\
\Gamma @>>\epsilon > A @>>\eta_{L} > \Gamma 
\end{CD}
\]
and a similar diagram involving $\chi \otimes 1$ and $\eta_{R}$.  This
is a great deal simpler than the corresponding diagram
in~\cite[Definition~A1.1.1(f)]{ravenel}. 

We recall that a \emph{$\Gamma$-comodule} is a left $A$-module $M$
together with a counital and coassociative coaction map $\psi \mathcolon
M\xrightarrow{}\Gamma \otimes_{A}M$ of left $A$-modules, where again
$\Gamma$ is an $A$-bimodule.  A map of comodules is a map of $A$-modules
that preserves the coaction, so we get a category $\Gamma \comod$ of
$\Gamma$-comodules.  

We then have the following
proposition~\cite[Proposition~2.2.8]{ravenel}, which explains the
importance of Hopf algebroids and comodules in algebraic topology.

\begin{proposition}\label{prop-homology-comodules}
Suppose $E$ is a ring spectrum such that $E_{*}E$ is a commutative ring
that is flat over $E_{*}$.  Then $(E_{*}, E_{*}E)$ is a Hopf algebroid,
and $E_{*}X$ is naturally an $E_{*}E$-comodule for $X$ a spectrum.  
\end{proposition}

Other good examples of Hopf algebroids include Hopf algebras, which are
just Hopf algebroids where $\eta_{L}=\eta_{R}$.  In particular, a
commutative ring $A$ can be thought of as the \emph{discrete} Hopf
algebroid $(A,A)$; a comodule over a discrete Hopf algebroid is just an
$A$-module.  Also, if $G$ acts on a commutative ring $R$ by ring
automorphisms, the ring of $R$-valued functions on $G$ is a Hopf
algebroid.  The right unit is defined by $\eta_{R}(r)(g)=g(r)$.  This is
dual to the twisted group ring $R[G]$.

We will summarize the properties of $\Gamma \comod$ in the following
theorem, but for it to make sense we need to recall some definitions.

\begin{definition}\label{defn-basics}
\begin{enumerate}
\item [(a)] A Hopf algebroid $\hopfalg$ is \textbf{flat} if $\eta_{R}$
makes $\Gamma$ into a flat $A$-module.  The conjugation shows that it is
equivalent to assume $\eta_{L}$ is flat.  Since $\epsilon$ is a right
inverse for both $\eta_{L}$ and $\eta_{R}$, both of these maps are then
faithfully flat.
\item [(b)] A category is \textbf{complete} if it has all small limits,
and \textbf{cocomplete} if it has all small colimits.  It is
\textbf{bicomplete} if it is both complete and cocomplete.  
\item [(c)] Given a regular cardinal $\lambda$, a category $\cat{I}$ is
said to be \textbf{$\lambda$-filtered} if every subcategory $\cat{J}$ of
$\cat{I}$ with fewer than $\lambda$ morphisms has an upper bound in
$\cat{I}$; that is, there is an object $C$ in $\cat{I}$ and a natural
transformation from the inclusion functor $\cat{J}\xrightarrow{}\cat{I}$
to the constant $\cat{J}$-diagram at $C$.  An object $M$ of a
cocomplete category $\cat{C}$ is said to be \textbf{$\lambda$-presented}
if $\cat{C}(M.-)$ commutes with $\lambda$-filtered colimits.
\item [(d)] Suppose $\cat{C}$ is a closed symmetric monoidal category
with monoidal structure $X\smash Y$, unit $A$, and closed structure
$F(X,Y)$.  An object $M$ in $\cat{C}$ is said to be \textbf{dualizable}
if the natural map $F(M,A)\smash X\xrightarrow{}F(M,X)$ is an
isomorphism for all $X\in \cat{C}$.  An object $M$ is said to be
\textbf{invertible} if there is an $N$ and an isomorphism $M\smash
N\cong A$.
\end{enumerate}
\end{definition}

\begin{theorem}\label{thm-comod-properties}
Suppose $\hopfalg$ is a flat Hopf algebroid.  Then $\Gamma \comod$ is a
bicomplete, closed symmetric monoidal abelian category.  We also have\uc 
\begin{enumerate}
\item [(a)] Filtered colimits are exact.  
\item [(b)] Given a regular cardinal $\lambda$ and a comodule $M$, $M$
is $\lambda$-presented if and only if $M$ is $\lambda$-presented as an
$A$-module.
\item [(c)] For any comodule $M$, there is a cardinal $\lambda$ such
that $M$ is $\lambda $-presented.  
\item [(d)] A comodule $M$ is dualizable if and only if it is projective
and finitely generated over $A$. 
\item [(e)] A comodule $M$ is invertible under the symmetric monoidal
product if and only if it is invertible as an $A$-module.  
\end{enumerate}
\end{theorem}

We denote the symmetric monoidal structure by $M\smash N$ with unit $A$
and the closed structure by $F(M,N)$.  

This theorem is a summary of the results
of~\cite[Section~1]{hovey-comodule-homotopy}.  We will just discuss some
of the issues that arise.  First of all, left adjoints are generally
easy to construct, since the forgetful functor from $\Gamma \comod$ to
$A\Mod$ is itself a left adjoint (Its right adjoint is the extended
comodule functor discussed in the following paragraph).  Thus, one
generally forms the left adjoint in $A\Mod$ and notices that it has a
natural comodule structure.  This is true for colimits and for the
symmetric monoidal structure $M\smash N$.  This is defined to be
$M\otimes_{A}N$, the tensor product of \textbf{left} $A$-modules, with
the coaction given as the composite
\[
M\otimes_{A}N \xrightarrow{\psi \otimes \psi } (\Gamma \otimes_{A}M)
\otimes_{A} (\Gamma \otimes_{A}N) \xrightarrow{g} \Gamma \otimes_{A}
M\otimes_{A} N   
\]
where $g(x\otimes m\otimes y\otimes n)=xy\otimes m\otimes n$.  

The key to constructing right adjoints is the extended comodule functor
from $A\Mod$ to $\Gamma \comod$ that takes $M$ to $\Gamma \otimes_{A}M$,
with coaction $\Delta \otimes 1$.  This is the right adjoint to the
forgetful functor.  As such, it is generally easy to define a desired
right adjoint $R$ on extended comodules.  For example, one can easily see
that we must define the product of extended comodules by 
\[
\prod^{\Gamma}(\Gamma \otimes_{A}M_{i}) \cong \Gamma \otimes_{A} \prod
M_{i}, 
\]
and the closed structure with target an extended comodule by 
\[
F(M, \Gamma \otimes_{A}N) \cong \Gamma \otimes_{A} \Hom_{A}(M,N).  
\]
It is less easy to see how one defines these right adjoints on maps
between extended comodules that are not necessarily extended maps, but
this can generally be done.  Having done this, we use the exact sequence
of comodules 
\[
0 \xrightarrow{} M \xrightarrow{\psi}\Gamma \otimes_{A}M
\xrightarrow{\psi g} \Gamma \otimes_{A} N,
\]
where $g\mathcolon \Gamma \otimes_{A}M\xrightarrow{}N$ is the cokernel
of $\psi$, to define $RM=\ker R(\psi g)$.  Since $R$ is supposed to be a
right adjoint, it must be left exact, so we must define $R$ in this
way.  

Note that if $M$ is an $A$-module and $N$ is a $\Gamma$-comodule, we
have the two tensor products $\Gamma \otimes_{A}(M\otimes_{A}N)$ and
$(\Gamma \otimes_{A}M)\smash N$.  It is useful to know that these are
the same~\cite{hovey-comodule-homotopy}.  

\begin{lemma}\label{lem-tensor}
Suppose $\hopfalg$ is a flat Hopf algebroid, $M$ is an $A$-module and
$N$ is a $\Gamma$-comodule.  Then there is a natural isomorphism of
comodules 
\[
(\Gamma \otimes_{A}M)\smash N \xrightarrow{} \Gamma
\otimes_{A}(M\otimes_{A}N).  
\]
\end{lemma}

Although Theorem~\ref{thm-comod-properties} indicates that the category
of $\Gamma$-comodules is a very well-behaved abelian category, one
obvious property is missing, and that is the existence of a set of
generators.  Recall that a set of objects $\cat{G}$ is said to
\textbf{generate} an abelian category $\cat{C}$ if, whenever $f$ is a
nonzero map in $\cat{C}$, there exists an object $G\in \cat{G}$ such
that $\cat{C}(G,f)$ is also nonzero.  For example, $A$ is a generator of
$A\Mod $.  This issue of generators is already complicated for Hopf
algebras; for a finite group $G$ and a field $k$, the natural generators
for the category of $k[G]$-modules (which is isomorphic to the category
of comodules over the ring of $k$-valued functions on $G$) are the
simple $k[G]$-modules.  There is no canonical description of these in
general.  However, any simple $k[G]$-module is finitely generated, and
of course projective, over $k$.  Referring to part~(d) of
Theorem~\ref{thm-comod-properties}, one might then expect that the set
of isomorphism classes of dualizable comodules forms a set of generators
for $\Gamma \comod$.  Sadly, this appears to be false in general, so we
need a hypothesis.  

\begin{definition}\label{defn-Adams}
A Hopf algebroid $\hopfalg$ is called an \textbf{Adams Hopf algebroid}
when $\Gamma$ is a filtered colimit of dualizable comodules.  
\end{definition}

This hypothesis is really due to
Adams~\cite[Section~III.13]{adams-blue}, who used it for the Hopf
algebroid $(E_{*}, E_{*}E)$ to set up universal coefficient spectral
sequences.  We learned it from~\cite{goerss-hopkins-comodules}, as well
as the following lemma.

\begin{lemma}\label{lem-Adams}
Suppose $\hopfalg$ is an Adams Hopf algebroid.  Then it is flat, and the
dualizable comodules generate the category of $\Gamma$-comodules.  
\end{lemma}

In categorical language, the category of comodules over an Adams Hopf
algebroid is a locally finitely presentable Grothendieck category.  

All of the Hopf algebroids that commonly arise in algebraic topology, as
well as all Hopf algebras over fields, are known to be
Adams~\cite[Section~1.4]{hovey-comodule-homotopy}.  

Note that it may be that one can take smaller sets of generators than
all of the dualizable comodules.  For example, the set $\{BP_{*}X_{n}
\}$ will serve as a set of generators for $BP_{*}BP$-comodules, where
$X_{n}$ is the $2n$-skeleton of $BP$.  However, $BP_{*}$ by itself is
definitely not a generator for the category of $BP_{*}BP$-comodules.  To
see this, let $PM$ be the set of primitives in a $BP_{*}BP$-comodule
$M$.  Then one can easily check that if $BP_{*}$ is a generator of the
category of $BP_{*}BP$-comodules then any comodule map that is
surjective after applying $P$ is in fact surjective.  In particular, the
map
\[
\bigoplus_{x\in PM} \Sigma^{|x|} BP_{*} \xrightarrow{f} M
\]
would be surjective.  This is easily seen to be false for
$M=BP_{*}(\mathbb{C}P^{2})$, for example.  

\begin{remark}
Theorem~\ref{thm-comod-properties} lists the good points of the category
$\Gamma \comod$; we now list some of the bad points of $\Gamma \comod$.

\begin{itemize}
\item [(a)] The forgetful functor from $\Gamma \comod$ to $A$-modules,
or even down to abelian groups, does NOT have a left adjoint; there is
no free comodule functor.  
\item [(b)] $\Gamma \comod$ does not, in general, have enough
projectives.  If we take $\hopfalg =(\Fp ,\cat{A})$, where $\cat{A}$
denotes the dual Steenrod algebra, it is generally believed that there
are no nonzero projective comodules.  
\item [(c)] If $\hopfalg$ is not a Hopf algebra, that is if $\eta_{L}$
and $\eta_{R}$ are not equal, then the forgetful functor from $\Gamma
\comod$ to $A\Mod$ is not in general surjective on objects.  For
example, there is no $BP_{*}BP$-comodule structure on $v_{n}^{-1}BP_{*}$
for $n>0$~\cite[Proposition~2.9]{johnson-yosimura}.  
\item [(d)] Products are not in general exact.  Hence the
inverse limit functor $\lim_{\Gamma}$ on sequences 
\[
\dotsb \xrightarrow{}M_{n}\xrightarrow{}\dotsb \xrightarrow{} M_{1}
\xrightarrow{} M_{0}
\]
may have nonzero derived functors $\lim_{\Gamma}^{i}$ for all $i>0$.  
\end{itemize}
\end{remark}

Because there are not enough projective comodules, the homological
algebra of comodules always involves injectives, or, better, relative
injectives.  Because the forgetful functor is exact, if $I$ is an
injective $A$-module, then $\Gamma \otimes_{A}I$ is an injective
$\Gamma$-comodule.  From this it is easy to check that there are enough
injectives.  However, injective $A$-modules are complicated, whereas
relative injectives have much better properties.  

A comodule $I$ is defined to be \textbf{relatively injective} if $\Gamma
\comod (-,I)$ takes $A$-split short exact sequences of comodules to
short exact sequences.  The following proposition sums up the properties
of relative injectives and is well-known; details can be found
in~\cite[Section~3.1]{hovey-comodule-homotopy}.  

\begin{proposition}\label{prop-rel-inj}
Suppose $\hopfalg$ is a flat Hopf algebroid.  
\begin{enumerate}
\item [(a)] The relatively injective $\Gamma$-comodules are the retracts
of extended comodules. 
\item [(b)] The coaction $\psi \mathcolon M\xrightarrow{}\Gamma
\otimes_{A}M$ defines an $A$-split embedding of $M$ into a relatively
injective comodule. 
\item [(c)] Relatively injective comodules are closed under coproducts
and products.  
\item [(d)] If $I$ is relatively injective, so is $I\smash M$ and
$F(M,I)$ for all comodules $I$.
\item [(e)] If $P$ is a comodule that is projective over $A$, and $I$ is
relatively injective, then $\Ext^{n}_{\Gamma}(P,I)=0$ for all $n>0$.  
\end{enumerate}
\end{proposition}

We take this proposition to mean that, to understand
$\Ext_{\Gamma}^{*}(M,N)$, we must simultaneously resolve $M$ by
comodules that are projective over $A$ and $N$ by relative injectives.
We return to this point in Section~\ref{sec-derived}.  

We close this section with a brief description of naturality. There is,
of course, a natural notion of a map $\Phi \mathcolon \hopfalg
\xrightarrow{}(B, \Sigma)$ of Hopf algebroids.  The map $\Phi$
corepresents a natural functor of groupoids, so consists of ring maps
$\Phi_{0}\mathcolon A\xrightarrow{}B$ and $\Phi_{1}\mathcolon \Gamma
\xrightarrow{}\Sigma$ satisfying certain conditions.  Such a map induces
a symmetric monoidal functor $\Phi_{*}\mathcolon \Gamma \comod
\xrightarrow{} \Sigma \comod$ that takes $M$ to $B\otimes_{A}M$, with
comodule structure given by the composite
\[
B\otimes_{A} M \xrightarrow{1\otimes \psi} B\otimes_{A}\Gamma
\otimes_{A}M \xrightarrow{g\otimes 1} \Sigma \otimes_{A} M \cong \Sigma
\otimes_{B} (B\otimes_{A}M)
\] 
where $g(b\otimes x)=b\Phi_{1}(x)$.  It is clear that $\Phi_{*}$
preserves colimits, so should have a right adjoint $\Phi^{*}$.  It does,
but, as usual, $\Phi^{*}$ is hard to define.  We define $\Phi^{*}(\Sigma
\otimes_{B}M)=\Gamma \otimes_{A}M$, and then extend this definition to
all $\Sigma$-comodules in the same way we did for the product of
comodules.  

An important new feature that arises in the study of Hopf algebroids is
the notion of weak equivalence.

\begin{definition}\label{defn-weak-equiv}
A map $\Phi$ of Hopf algebroids is defined to be a \textbf{weak
equivalence} if $\Phi_{*}$ is an equivalence of categories.  
\end{definition}

If $\Phi$ is a weak equivalence between discrete Hopf algebroids, then
$\Phi$ is an isomorphism, but a central point of the author's work on
Hopf algebroids is that there are important non-trivial weak
equivalences of Hopf algebroids that are not discrete.  

In general, we have the following characterization of weak
equivalences. 

\begin{theorem}\label{thm-weak-equiv}
The map $\Phi$ of Hopf algebroids is a weak equivalence if and only if
the composite 
\[
A \xrightarrow{\eta_{R}} \Gamma  \xrightarrow{\Phi_{0}\otimes 1}
B\otimes_{A} \Gamma  
\]
is a faithfully flat ring extension and the map 
\[
B\otimes_{A} \Gamma \otimes_{A} B \xrightarrow{} \Sigma 
\]
that takes $b\otimes x\otimes b'$ to
$\eta_{L}(b)\Phi_{1}(x)\eta_{R}(b')$ is a ring isomorphism.  
\end{theorem}

The ``if''half of this theorem is the main result of~\cite{hovey-hopf}.
The ``only if'' half is much easier and was proven
in~\cite{hovey-strickland-comodules}.  

This theorem has a better formulation.  Hollander~\cite{hollander} has
constructed a model structure on presheaves of groupoids on a
Grothendieck topology $\cat{C}$; the fibrant objects are stacks.  In
particular, we can take our Grothendieck site to be the flat topology on
$\Aff$, the opposite category of commutative rings (with a cardinality
bound so we get a small category).  In this topology, a cover of $R$ is
a finite collection of flat extensions $S_{i}$ of $R$ such that $\prod
S_{i}$ is faithfully flat over $R$.  Any Hopf algebroid $\hopfalg$
defines a presheaf of groupoids $\Spec \hopfalg$ by definition; this
presheaf is in fact a sheaf in the flat topology by faithfully flat
descent~\cite{hovey-hopf}.  We can then rephrase
Theorem~\ref{thm-weak-equiv} as follows.

\begin{corollary}\label{cor-weak-equiv}
A map $\Phi$ of Hopf algebroids is a weak equivalence if and only if
$\Spec \Phi$ is a weak equivalence of sheaves of groupoids in the flat
topology.  
\end{corollary}

This point of view suggests that one should reconsider the results of
this section for quasi-coherent sheaves over a sheaf of groupoids, since
a quasi-coherent sheaf over $\Spec \hopfalg$ is the same thing as a
$\Gamma$-comodule~\cite{hovey-hopf}.  We have not carried out this
program.  One reason for this is that we don't see any clear
applications.  But another reason is that we do not know whether the
Adams condition is invariant under weak equivalence.  The difficulty is
that, while dualizable comodules and filtered colimits are preserved by
weak equivalences, $\Gamma$ is not.  One could simply demand that
dualizable sheaves generate the category of quasi-coherent sheaves, but
again we do not know if this is sufficient to provide a useful theory,
or even whether it holds for interesting non-affine sheaves of
groupoids.

It would be interesting to know if there are equivalences of categories
of comodules that are not given by maps of Hopf algebroids, as occurs in
Morita theory.  Since Hopf algebroids are a generalization of
\textbf{commutative} rings, and there are no non-trivial Morita
equivalences of commutative rings, it is reasonable to guess that every
equivalence of categories of comodules is a zig-zag of weak
equivalences.  

\section{Landweber exact algebras}\label{sec-Landweber}

The object of this section is to study the relation between
$\Gamma$-comodules and $\Gamma_{B}$-comodules, where $B$ is a Landweber
exact $A$-algebra.  The main application is to the relation between
$BP_{*}BP$-comodules and $E(n)_{*}E(n)$-comodules.  In particular, we
sketch the proof of Theorem~\ref{main-A} and its corollaries in this
section.  More details can be found
in~\cite{hovey-strickland-comodules}.  

Given a Hopf algebroid $\hopfalg$, an $A$-algebra $B$ is said to be
\textbf{Landweber exact} over $A$ if $B\otimes_{A}(-)$ takes exact
sequences of $\Gamma$-comodules to exact sequences of $B$-modules.  This
is called Landweber exactness because Landweber gave a characterization
of Landweber exact $BP_{*}$-algebras in his famous Landweber exact
functor theorem~\cite{land-exact}.  One can check that $B$ is Landweber
exact over $A$ if and only if the composite
\[
A \xrightarrow{\eta_{R}} \Gamma  \xrightarrow{} B\otimes_{A} \Gamma 
\]
is a flat ring extension.  This condition is reminiscent of the
characterization of weak equivalences given in
Theorem~\ref{thm-weak-equiv}.  We can make it even more so by defining 
\[
\Gamma_{B}=B \otimes_{A}\Gamma \otimes_{A}B.  
\]

We then have the following lemma, which is easy to prove but can also be
found in~\cite{hovey-strickland-comodules}.  

\begin{lemma}\label{lem-Landweber}
Suppose $\hopfalg$ is a Hopf algebroid and $B$ is an $A$-algebra.  Then
$(B, \Gamma_{B})$ is a Hopf algebroid, and the evident map $\hopfalg
\xrightarrow{}(B, \Gamma_{B})$ is a map of Hopf algebroids.  If
$\hopfalg$ is flat and $B$ is Landweber exact, then $(B, \Gamma_{B})$ is
a flat Hopf algebroid.
\end{lemma}

Thus, if $B$ is Landweber exact over $A$, the map 
\[
\Phi \mathcolon \hopfalg \xrightarrow{}(B, \Gamma_{B})
\]
is almost a weak equivalence, in that 
\[
A \xrightarrow{\eta_{R}} \Gamma  \xrightarrow{} B\otimes_{A} \Gamma 
\]
is flat, and 
\[
B\otimes_{A} \Gamma \otimes_{A} B \xrightarrow{} \Gamma_{B}
\]
is an isomorphism.  The only thing stopping $\Phi$ from being a weak
equivalence is that $B\otimes_{A}(-)$ may not be faithful on the
category of $\Gamma$-comodules.  The idea of the following theorem,
proved in~\cite{hovey-strickland-comodules}, is that we can force $\Phi$
to be faithful by localizing the category of $\Gamma$-comodules.  

\begin{theorem}\label{thm-Giraud}
Suppose $\hopfalg$ is a flat Hopf algebroid, and $B$ is Landweber exact
over $A$.  Then the map 
\[
\Phi \mathcolon \hopfalg \xrightarrow{} (B, \Gamma_{B})
\]
yields an equivalence 
\[
\Phi_{*} \mathcolon L_{\cat{T}}(\Gamma \comod)
\xrightarrow{}\Gamma_{B}\comod 
\]
where $L_{\cat{T}}(\Gamma \comod)$ is the localization of $\Gamma
\comod$ with respect to the hereditary torsion theory $\cat{T}$
consisting of all $\Gamma$-comodules $M$ such that $B\otimes_{A}M=0$.  
\end{theorem}

A \textbf{hereditary torsion theory} is just a full subcategory closed under
subobjects, quotient objects, extensions, and arbitrary direct sums.
The localization $L_{\cat{T}}$ is obtained by inverting all maps $f$ of
$\Gamma$-comodules whose kernel and cokernel are in $\cat{T}$.  Note
that the Hopf algebroids that arise in algebraic topology are graded, so
$B$ will be a graded $A$-algebra, and our hereditary torsion theories
will also be graded, in the sense that $M$ is in $\cat{T}$ if and only
if all shifts of $M$ are in $\cat{T}$.  

Because it is so surprisingly easy, we will give the proof of
Theorem~\ref{thm-Giraud}.

\begin{proof}
Consider the natural transformation 
\[
\epsilon_{M}\mathcolon \Phi_{*}\Phi^{*}M\xrightarrow{}M.
\]
We claim that this map is a natural isomorphism.  One can check this by
calculation for extended $\Sigma$-comodules $M$.  Since $\epsilon$ is a
natural transformation of left exact functors (because $B$ is Landweber
exact), and every $\Sigma$-comodule is the kernel of a map of extended
comodules, $\epsilon_{M}$ is an isomorphism for all $M$.  

After this, the rest of the proof of Theorem~\ref{thm-Giraud} is purely
formal.  A priori, the category $L_{\cat{T}}(\Gamma \comod)$ may not be
an actual category, since it may not have small Hom sets, always a
danger with localization.  However, it does exist in a higher universe.
The natural transformation
\[
\eta_{M} \mathcolon M\xrightarrow{}\Phi^{*}\Phi_{*}M
\]
becomes an isomorphism upon applying $\Phi_{*}$, and therefore, since
$\Phi_{*}$ is exact, the kernel and cokernel of $\eta_{M}$ are in
$\cat{T}$.  This gives us the desired equivalence, and incidentally
shows that $L_{\cat{T}}(\Gamma \comod)$ actually does exist as an honest
category, since it is equivalent to $\Gamma_{B}\comod$.  
\end{proof}

Theorem~\ref{thm-Giraud} gives the following corollary, which
it is difficult to imagine proving directly.  

\begin{corollary}\label{cor-primitive}
Suppose $\hopfalg$ is a flat Hopf algebroid, $B$ is Landweber exact over
$A$, and every nonzero $\Gamma $-comodule has a primitive.  Then every
nonzero $\Gamma_{B}$-comodule has a primitive.  In particular, every
$E(n)_{*}E(n)$-comodule has a primitive. 
\end{corollary}

This corollary is immediate, as $L_{\cat{T}}(\Gamma \comod)$ is the
full subcategory of $\Gamma \comod$ consisting of the local objects.  

To get further information, we need to identify the hereditary torsion
theories that can arise.  Let $\cat{T}_{n}$ denote the collection of all
$v_{n}$-torsion $BP_{*}BP$-comodules, so that $\cat{T}_{0}$ is the
collection of all $p$-torsion comodules, and $\cat{T}_{-1}$ is the
collection of all comodules.  One can easily check that $\cat{T}_{n}$ is
a hereditary torsion theory.  It is less obvious, but true, that
$\cat{T}_{n}$ is the smallest graded hereditary torsion theory
containing $BP_{*}/I_{n+1}$.

\begin{theorem}\label{thm-torsion}
Let $\cat{T}$ be a graded hereditary torsion theory of
$BP_{*}BP$-comodules.  If $\cat{T}$ contains a nonzero finitely presented
comodule, then $\cat{T}=\cat{T}_{n}$ for some $n$.
\end{theorem}

This theorem is proved in~\cite{hovey-strickland-comodules}, using the
ideas behind the Landweber filtration theorem.  Note that this theorem
explains why $L_{n}$ and $L_{n}^{f}$ agree on the category of
$BP_{*}BP$-comodules.  The only reasonable definition of $L_{n}$ is
localization with respect to the hereditary torsion theory of all
comodules $M$ such that $E(n)_{*}\otimes_{BP_{*}}M=0$, and $L_{n}^{f}$
is localization with respect to the hereditary torsion theory generated
by $BP_{*}/I_{n+1}$.  Both of these torsion theories are $\cat{T}_{n}$.  

Because of this theorem, it is natural to make the following definition.

\begin{definition}\label{defn-height}
Suppose $B$ is a $BP_{*}$-algebra.  Define the \textbf{height} of $B$ to
be the largest integer $n$ such that $B/I_{n}$ is nonzero.  If there is
no such $n$, define the height of $B$ to be infinite.  
\end{definition}

From a formal group law point of view, the height of $B$ is the largest
possible height of any specialization of the formal group law of $B$.
So the height of $E(n)$ is $n$, and the height of $BP$ itself is
$\infty$.  

Here is the main theorem of~\cite{hovey-strickland-comodules}.

\begin{theorem}\label{thm-equiva}
Let $(A, \Gamma)=(BP_{*},BP_{*}BP)$, and suppose $B$ is a graded
Landweber exact $A$-algebra of height $n\leq \infty $.  Then the
functor $M\mapsto B\otimes_{A}M$ defines an equivalence of
categories
\[
L_{n}(\Gamma \comod ) \xrightarrow{} \Gamma_{B}\comod,
\]
where $L_{n}$ is localization with respect to $\cat{T}_{n}$ for
$n<\infty$ and $L_{\infty}$ is the identity localization.  
\end{theorem}

This theorem is almost a corollary of Theorem~\ref{thm-Giraud} and
Theorem~\ref{thm-torsion}, except for the infinite height case.  We do
not have a classification of graded hereditary torsion theories of
$BP_{*}BP$-comodules that do not contain a nonzero finitely presented
comodule.  There are probably uncountably many such torsion theories.
However, if $B$ is Landweber exact and $B\otimes_{A}M=0$ for some
nonzero $M$, then $B/I_{\infty}B=0$.  Indeed, $M$ must have a nonzero
primitive, and so we conclude by Landweber exactness that $B/IB=0$ for
some proper invariant ideal $I$.  Since $I\subseteq I_{\infty}$, it
follows that $B/I_{\infty}B=0$.  But this means that $1\in I_{\infty}B$,
so $1\in I_{n}B$ for some $n$.  Thus $B\otimes_{A}A/I_{n}=0$.  This
proves that if $B$ has infinite height, then $B\otimes_{A}(-)$ does not
kill any nonzero $BP_{*}BP$-comodules.   

The following corollary is immediate. 

\begin{corollary}\label{cor-equiv}
Let $\hopfalg =(BP_{*},BP_{*}BP)$, and suppose $B$ and $B'$ are both
Landweber exact $A$-algebras of the same height.  Then the category of
$\Gamma_{B}$-comodules is equivalent to the category of
$\Gamma_{B'}$-comodules.  
\end{corollary}

In particular, the category of $E(n)_{*}E(n)$-comodules is equivalent to
the category of $v_{n}^{-1}BP_{*}(v_{n}^{-1}BP)$-comodules, even though
the category of $E(n)_{*}$-modules is very different from the category
of $v_{n}^{-1}BP_{*}$-modules.  One way to think of the Miller-Ravenel
change of rings theorem~\cite[Theorem~3.10]{miller-ravenel} as an
isomorphism of certain $\Ext$ groups in these two categories.  This is
now obvious; the $\Ext$ groups are isomorphic because the categories
they are taken in are equivalent.

We point out that the equivalence of categories of comodules in
Corollary~\ref{cor-equiv} is in fact induced by a zig-zag of weak
equivalences of Hopf algebroids.  Any map $B\xrightarrow{}B'$ of
Landweber exact $BP_{*}$-algebras of the same height induces a weak
equivalence of Hopf algebroids 
\[
(B, \Gamma_{B})\xrightarrow{}(B', \Gamma_{B'}),  
\]
where we are still denoting $BP_{*}BP$ by $\Gamma$.  If $B$ and $B'$ are
Landweber exact $BP_{*}$-algebras of the same height, there may not be a
map of $BP_{*}$-algebras between them.  However, if we let
$C=B\otimes_{BP_{*}}\Gamma \otimes_{BP_{*}}B'$, then $C$ has a left and
right $BP_{*}$-algebra structure, which we denote by $C_{L}$ and
$C_{R}$.  There are maps of $BP_{*}$-algebras $B\xrightarrow{}C_{L}$ and
$B'\xrightarrow{}C_{R}$, and conjugation induces an isomorphism
$C_{L}\xrightarrow{}C_{R}$.  Since $C_{L}$ is also Landweber exact, of
the same height as $B$ and $B'$, this yields the desired zig-zag of weak
equivalences.

To further understand the structure of the category of
$E(n)_{*}E(n)$-comodules, we would like to understand the localization
functor $L_{n}$ better.  The following theorem is a summary of the
results of~\cite{hovey-strickland-derived}, and is joint work of the author
and Strickland.  

\begin{theorem}\label{thm-localization}
\begin{enumerate}
\item [(a)] A comodule $M$ is $L_{n}$-local if and only if 
\[
\Hom_{BP_{*}BP}^{*}(BP_{*}/I_{n+1},M) =
\Ext^{1,*}_{BP_{*}BP}(BP_{*}/I_{n+1},M)=0,
\]
which is true if and only if 
\[
\Hom_{BP_{*}}^{*}(BP_{*}/I_{n+1},M) =
\Ext^{1,*}_{BP_{*}}(BP_{*}/I_{n+1},M)=0,
\]
\item [(b)] $L_{n}$, thought of as an endofunctor of the category of
$BP_{*}BP$-comodules, is left exact and preserves finite limits,
filtered colimits, and arbitrary direct sums.  It has right derived
functors which we denote $L_{n}^{i}$. 
\item [(c)] For $i>0$, $L_{n}^{i}(M)$ is isomorphic to the $i+1$st local
cohomology group of the $BP_{*}$-module $M$ with respect to $I_{n+1}$.
In particular, $L_{n}^{i}(M)=0$ for $i>n$.  
\item [(d)] Suppose $m<n$ and $M$ is a $v_{m-1}$-torsion comodule on
which $(v_{m}, v_{m+1})$ is a regular sequence.  Then $M$ is
$L_{n}$-local.
\item [(e)] If $v_{m}$ acts invertibly on a comodule $M$ for some $m\leq
n$, then $M$ is $L_{n}$-local and $L^{i}_{n}M=0$ for all $i>0$.  
\item [(f)] If a comodule $M$ is $v_{n-1}$-torsion, then
$L_{n}M=v_{n}^{-1}M$.
\item [(g)] We have 
\[
L_{n}(BP_{*}/I_{k})= \begin{cases}
BP_{*}/I_{k} & k <n, \\
v_{n}^{-1}BP_{*}/I_{n} & k=n, \\
0 & k>n.
\end{cases}
\]
and, for $i,n>0$, 
\[
L_{n}^{i}(BP_{*}/I_{k}) = \begin{cases}
BP_{*}/(p, v_{1},\dotsc ,v_{k-1}, v_{k}^{\infty},\dotsc ,v_{n}^{\infty}) & i=n-k>0, \\
0 & \text{otherwise}.
\end{cases}
\]
\end{enumerate}
\end{theorem}

These derived functors $L_{n}^{i}$ can be used to compute
$BP_{*}(L_{n}X)$ from $BP_{*}X$ by means of a spectral sequence.  

\begin{theorem}\label{thm-spectral}
Let $X$ be a spectrum.  There is a natural spectral sequence
$E_{*}^{**}(X)$ with $d_{r}\mathcolon
E_{r}^{s,t}\xrightarrow{}E_{r}^{s+t,t+r-1}$ and $E_{2}$-term
$E_{2}^{s,t}(X)\cong (L^{s}_{n}BP_{*}X)_{t}$, converging to
$BP_{t-s}(L_{n}X)$.  This is a spectral sequence of
$BP_{*}BP$-comodules, in the sense that $E_{r}^{s,*}$ is a
graded $BP_{*}BP$-comodule for all $r\geq 2$ and $d_{r}\mathcolon
E_{r}^{s,*}\xrightarrow{}E_{r}^{s+r,*}$ is a $BP_{*}BP$-comodule map of
degree $r-1$.  Furthermore, every element in $E_{2}^{0,*}$ that comes
from $BP_{*}X$ is a permanent cycle.  
\end{theorem}

This theorem is proved in~\cite{hovey-strickland-derived}.  It is
very closely related to the local cohomology spectral sequence of
Greenlees~\cite{greenlees-spectral} and Greenlees
and May~\cite{greenlees-may-completions}.  One way of putting it is that
we show that the Greenlees spectral sequence is a spectral sequence
of comodules in this case.  When $X=S^{0}$, this implies that the
spectral sequence collapses with no extensions, and we recover Ravenel's
computation~\cite{rav-loc} of $BP_{*}L_{n}S^{0}$.  

We now derive some corollaries of Theorem~\ref{main-B}, proved
in~\cite{hovey-strickland-comodules}.  

\begin{corollary}\label{cor-local-compute}
Let $\hopfalg =(BP_{*}, BP_{*}BP)$, and suppose $B$ is a Landweber exact
$A$-algebra of height $n$.  Then 
\[
\Gamma_{B}\comod (B, B/I_{m}) = \begin{cases}  \Z_{(p)} & n>m=0 \\
\Q & n=m=0 \\
\Fp [v_{m}] & n>m>0 \\
\Fp [v_{m}, v_{m}^{-1}] & n=m>0 .
\end{cases}
\]
\end{corollary}

\begin{proof}
Let $\Phi \mathcolon (A, \Gamma)\xrightarrow{} (B, \Gamma_{B})$ be the
evident map of Hopf algebroids, so that $L_{n}=\Phi^{*}\Phi_{*}$ by
Theorem~\ref{thm-equiva}.  Then we have
\[
\Gamma_{B}\comod (B, B/I_{m}) \cong \Gamma \comod (A, \Phi^{*}(B/I_{m}))
\cong \Gamma \comod (A, L_{n}(A/I_{m})), 
\]
so the result follows from Theorem~\ref{thm-localization} and the
analogous calculation for $B=BP_{*}$~\cite[Theorem~4.3.2]{ravenel}.  
\end{proof}

This corollary in turn gives rise to the expected structural results
about $\Gamma_{B}$-comodules, proved
in~\cite{hovey-strickland-comodules}.  

\begin{theorem}\label{thm-Landweber}
Let $(A, \Gamma)=(BP_{*}, BP_{*}BP)$, and suppose $B$ is a Landweber
exact $A$-algebra of height $n$.  
\begin{enumerate}
\item [(a)] If $I$ is an invariant radical ideal in $B$, then $I=I_{m}$
for some $m\leq n$.  
\item [(b)] If $M$ is a finitely presented $\Gamma_{B}$-comodule, then
there is a filtration 
\[
0=M_{0} \subseteq M_{1} \subseteq \dotsb \subseteq M_{t} =M
\]
of $M$ by subcomodules such that, for each $s\leq t$, there is a $k\leq
n$ and an $r$ such that $M_{s}/M_{s-1}\cong s^{r}B/I_{k}B$.  
\end{enumerate}
\end{theorem}

\begin{proof}
Note that the theorem is invariant under the equivalences of categories
of Theorem~\ref{thm-equiva}.  Thus we might as well assume that
$B=E(n)_{*}$ or $BP_{*}$, and in case $B=BP_{*}$ the result is due to
Landweber~\cite{land-exact}.  So we assume $B=E(n)_{*}$.  Now, for
part~(a), suppose $I$ is an invariant radical ideal, and choose the
largest $k$ such that $I_{k}\subseteq I$.  If $I_{k}\neq I$, the
comodule $I/I_{k}$ must have a nonzero primitive $y$.  This primitive
must also be a primitive in $E(n)_{*}/I_{k}$.  By
Corollary~\ref{cor-local-compute}, it must be a power of $v_{k}$.  Since
$I$ is radical, this means that $I_{k+1}\subseteq I$.  This
contradication implies that $I_{k}=I$.  

For part~(b), we construct the filtration $M_{i}$ by induction on $i$,
taking $M_{0}=0$.  Having built $M_{i}$, if $M_{i}\neq M$, we choose a
nonzero primitive $y$ in $M/M_{i}$.  We claim that some multiple $z$ of
$y$ is a primitive whose annihilator is $I_{k}$ for some $k\leq n$.
Indeed, if $y$ is $p$-torsion, then we can multiply $y$ by a power of
$p$ to obtain a nonzero primitive $y_{1}$ with $py_{1}=0$.  Since
$v_{1}$ is a primitive mod $p$, if $y_{1}$ is $v_{1}$-torsion we can
multiply $y_{1}$ by power of $v_{1}$ to obtain a nonzero primitive
$y_{2}$ with $py_{2}=v_{1}y_{2}=0$. Continuing in this fashion, we end
up with a primitive $z$ such that $I_{k}z=0$ and $z$ is not
$v_{k}$-torsion.  Corollary~\ref{cor-local-compute} implies that $\Ann
z=I_{k}$, as required.  

We now choose an element $w$ in $M$ whose image in $M/M_{i}$ is $z$, and
let $M_{i+1}$ denote the subcomodule generated by $M_{i}$ and $w$.  Then
$M_{i+1}/M_{i}\cong s^{r}E(n)_{*}/I_{k}$, where $r$ is the degree of
$z$.  Since $M$ is finitely presented over the Noetherian ring
$E(n)_{*}$, this process must stop, and so $M_{t}=M$ for some $t$.   
\end{proof}

\section{The stable homotopy category $\stable{\Gamma }$}\label{sec-derived}

The object of this section is to discuss the construction and basic
properties of the stable homotopy category $\stable{\Gamma}$ associated
to an Adams Hopf algebroid $\hopfalg$.  In practice, we are most
interested in $\stable{E_{*}E}$ for $E$ a commutative Landweber exact
ring spectrum.  This means our Hopf algebroids should be graded, and the
stable homotopy category $\stable{\Gamma}$ should be bigraded.  However,
the grading just adds unnecessary complexity to the notation, so we will
forget about it for most of this section.

\subsection*{Construction of $\stable{\Gamma }$}

The usual derived category $\cat{D}(\Gamma)$ of $\Gamma$-comodules is
obtained from $\Ch{\Gamma}$, the category of unbounded chain complexes
of $\Gamma$-comodules, by inverting the homology isomorphisms.  This is
not the right thing to do to form $\stable{\Gamma}$.  To see this, note
that there is a model structure on $\Ch{\Gamma}$ whose homotopy category
is the derived category, as there is on $\Ch{\cat{A}}$ for any
Grothendieck category $\cat{A}$~\cite{beke}.  The cofibrations in this
model structure are the monomorphisms, the fibrations are the
epimorphisms with DG-injective kernel, and the weak equivalences are the
homology isomorphisms.  Here a complex $X$ is \textbf{DG-injective} if
each $X_{n}$ is injective and every map from an exact complex to $X$ is
chain homotopic to $0$.  In particular, bounded above complexes of
injectives are DG-injective, from which it follows that
\[
\cat{D}(\Gamma)(A, A)_{*} = \Ext_{\Gamma}^{*}(A,A).
\]
Here we are thinking of $A$ as a complex concentrated in degree $0$.  It
is the analog of the sphere $S$, since it is the unit of $\smash$.  Now,
there are often non-nilpotent elements in $\Ext^{i}_{\Gamma}(A,A)$ for
$i>0$.  For example, the well-known element $\beta_{1}$ is non-nilpotent
in $\Ext^{*}_{BP_{*}BP}(BP_{*},BP_{*})$ for $p>2$.  Since $\beta_{1}$
corresponds to a self-map of $BP_{*}$ in $\cat{D}(BP_{*}BP)$, we can
form $\beta_{1}^{-1}BP_{*}$.  But this complex has trivial homology, so
is $0$ in $\cat{D}(BP_{*}BP)$ even though $\beta_{1}$ is not nilpotent.
This is not good for several reasons; we should be able to see
$\beta_{1}^{-1}BP_{*}$ because it is an important object, and the fact
that we can't also implies that $A$ is not a small object of
$\cat{D}(BP_{*}BP)$.

We should be inverting homotopy isomorphisms, not homology
isomorphisms.  To do this, we need to define the homotopy groups. 

\begin{definition}\label{defn-homotopy}
Suppose $\hopfalg$ is an Adams Hopf algebroid, $X\in \Ch{\Gamma}$, and
$P$ is a dualizable $\Gamma$-comodule.  We define the \textbf{homotopy
groups of $X$ with coefficients in $P$} by 
\[
\pi_{n}^{P}(X) = H_{-n}(\Gamma \comod (P, LA\smash X)),
\]
where $LA$ is the cobar resolution of $A$.  We define a map $f$ of
complexes to be a \textbf{homotopy isomorphism} if $\pi_{n}^{P}(f)$ is
an isomorphism for all $n$ and all dualizable comodules $P$.  The
\textbf{stable homotopy category of $\Gamma$}, $\stable{\Gamma}$ is
defined to be the category obtained from $\Ch{\Gamma}$ by inverting the
homotopy isomorphisms.  
\end{definition}

The reason for the sign in the definition of $\pi_{n}^{P}$ is so that 
\[
\pi_{n}^{P}(M)\cong \Ext^{n}_{\Gamma}(P,M)
\]
for a comodule $M$.  Homotopy groups and homotopy isomorphisms satisfy
the expected properties~\cite{hovey-comodule-homotopy}.

\begin{proposition}\label{prop-homotopy-groups}
Suppose $\hopfalg$ is an Adams Hopf algebroid, and $P$ is a dualizable
comodule.  
\begin{enumerate}
\item [(a)] A short exact sequence of complexes induces a long exact
sequence in the homotopy groups $\pi_{*}^{P}(-)$.  Note that the
boundary map raises dimension by one, because of the sign in the
definition of $\pi_{n}^{P}(-)$.
\item [(b)] Homotopy groups commute with filtered colimits, so homotopy
isomorphisms are closed under filtered colimits.  
\item [(c)] Every chain homotopy equivalence is a homotopy isomorphism,
and every homotopy isomorphism is a homology isomorphism.  
\item [(d)] The natural map $X\xrightarrow{}LA\smash X$ is a homotopy
isomorphism for all $X$.  
\end{enumerate}
\end{proposition}

To understand $\stable{\Gamma}$, and for that matter to even see that it
is a category at all, we need a model structure on $\Ch{\Gamma}$ in
which the weak equivalences are the homotopy isomorphisms.  This is the
main goal of~\cite{hovey-comodule-homotopy}, where the following theorem
is proved. 

\begin{theorem}\label{thm-model}
Suppose $\hopfalg$ is an Adams Hopf algebroid.  Then there is a proper
symmetric monoidal model structure on $\Ch{\Gamma}$ in which the weak
equivalences are the homotopy isomorphisms.  Furthermore, if $X$ is
cofibrant, then $X\smash (-)$ preserves homotopy isomorphisms.  
\end{theorem}

We call this model structure the \textbf{homotopy model structure}.  The
cofibrations in the homotopy model structure are dimensionwise split
monomorphisms with cofibrant cokernel.  If $X$ is cofibrant, then each
$X_{n}$ is projective over $A$.  More precisely, if $X$ is cofibrant,
then $X$ is a retract of a complex $Y$ that admits a filtration $Y^{i}$
such that each map $Y^{i}\xrightarrow{}Y^{i+1}$ is a dimensionwise split
monomorphism and the quotient $Y^{i+1}/Y^{i}$ is a complex of relatively
projective comodules with trivial differential.  Here a comodule is
\textbf{relatively projective} if it is a retract of a direct sum of
dualizable comodules.

A characterization of the fibrations is given
in~\cite{hovey-comodule-homotopy}.  Fibrations are of course
surjective.  Every complex of relative injectives
is fibrant, and every fibrant complex is equivalent in a precise sense
to a complex of relative injectives.  A fibrant replacement of $X$ is
given by $LB\smash X$.  

The following theorem is also proved in~\cite{hovey-comodule-homotopy}.

\begin{theorem}\label{thm-model-natural}
The homotopy model structure is natural, in the sense that a map $\Phi
\mathcolon \hopfalg \xrightarrow{} (B, \Sigma)$ induces a left Quillen
functor $\Phi_{*} \mathcolon \Ch{\Gamma}\xrightarrow{}\Ch{\Sigma}$ of the
homotopy model structures.  Furthermore, if $\Phi$ is a weak equivalence,
then $\Phi_{*}$ is a strong Quillen equivalence, in the sense that both
$\Phi_{*}$ and $\Phi^{*}$ preserve and reflect homotopy isomorphisms.  
\end{theorem}

\subsection*{Global properties of $\stable{\Gamma}$}

We now establish some of the essential properties of $\stable{\Gamma}$.
At this point, we begin to use some of the standard notational
conventions of ordinary stable homotopy.  Thus, we will begin using $S$
for the image of $A$ in $\stable{\Gamma}$, thinking of it as analogous
to the usual zero-sphere.  Similarly, we will sometimes use $[X,Y]_{*}$
for graded maps in $\stable{\Gamma}$.  In practice, $\Gamma$ is usually
graded and so $\stable{\Gamma}$ is bigraded.  However, the internal
suspension in the category of $\Gamma$-comodules is usually not
relevant, so we tend to omit it from the notation.  

\begin{theorem}\label{thm-model-homotopy}
Suppose $\hopfalg$ is an Adams Hopf algebroid.  The category
$\stable{\Gamma }$ is a closed symmetric monoidal triangulated
category.  The dualizable comodules form a set of small, dualizable,
weak generators for $\stable{\Gamma}$.  
\end{theorem}

This theorem is really a corollary of Theorem~\ref{thm-model} and
general facts about model categories.  It is proved
in~\cite{hovey-comodule-homotopy}; another way to say it is that
$\stable{\Gamma}$ is a unital algebraic stable homotopy category in the
sense of~\cite{hovey-axiomatic}.

One drawback of $\stable{\Gamma}$ is that it is not in general
monogenic.  That is, $A$ and its suspensions are not generally enough to
generate the whole category.  This is unavoidable even for Hopf
algebras.  Indeed, if $G$ is a finite group and $k$ is a field, the
stable homotopy category of the Hopf algebra of functions from $G$ to
$k$ is closely related to the stable module category much studied in
modular representation theory~\cite{benson-book-1}, as explained
in~\cite{hovey-axiomatic}.  If $G$ is a $p$-group, the stable module
category is monogenic, but not in general.  

However, one certainly expects $\stable{BP_{*}BP}$ and
$\stable{E(n)_{*}E(n)}$ to be monogenic, so we need a condition on
$\hopfalg$ that will ensure that $\stable{\Gamma}$ is monogenic.  Recall
that a full subcategory $\cat{D}$ of an abelian category is called
\textbf{thick} if it is closed under retracts and, whenever two out of
three terms in a short exact sequence are in $\cat{D}$, so is the
third. 

\begin{proposition}\label{prop-monogenic}
Suppose $\hopfalg$ is an Adams Hopf algebroid, and every dualizable
comodule is in the thick subcategory generated by $A$.  Then
$\stable{\Gamma}$ is monogenic.  
\end{proposition}

This proposition is proved in~\cite{hovey-comodule-homotopy}.  To apply
it, we note that the filtration theorem for $E_{*}E$-comodules, part~(b)
of Theorem~\ref{thm-Landweber}, implies that every finitely presented
$E_{*}E$-comodule is in the thick subcategory generated by $E_{*}$ when
$E$ is a commutative ring spectrum that is Landweber exact over
$BP_{*}$.  Thus we get the following corollary.

\begin{corollary}\label{cor-monogenic}
Suppose $E$ is a commutative ring spectrum that is Landweber exact over
$BP$.  Then $\stable{E_{*}E}$ is monogenic in the bigraded sense.  In
particular, a map $f$ in $\Ch{E_{*}E}$ is a homotopy isomorphism if and
only if $\pi_{n,k}(f)= \pi_{n}^{s^{k}E_{*}}(f)$ is an isomorphism for
all $n$ and $k$.
\end{corollary}

The bigrading arises because we can suspend a complex $X$ either
internally, by suspending each graded comodule $X_{n}$, or externally by
suspending the complex $X$.  

We would like to understand the relation between a comodule $M$ and its
image in $\stable{\Gamma}$.  The following proposition is proved
in~\cite{hovey-comodule-homotopy}.  

\begin{proposition}\label{prop-comodule-stable}
Suppose $\hopfalg$ is an Adams Hopf algebroid.  
\begin{enumerate}
\item [(a)] A short exact sequence of comodules, or even complexes,
gives rise to a cofiber sequence in $\stable{\Gamma}$.
\item [(b)] If $M$ is in the thick subcategory generated by $A$, then
$M$ is a small object of $\stable{\Gamma}$.  
\item [(c)] If $M$ and $N$ are comodules, then there is a natural map 
\[
\Ext^{k}_{\Gamma}(M,N) \xrightarrow{} \stable{\Gamma}(M,N)^{k}
\]
that is an isomorphism for $M$ in the thick subcategory generated by
$A$.  
\end{enumerate}
\end{proposition}

Part~(b) is not actually proved in~\cite{hovey-comodule-homotopy}, but
follows immediately from part~(a).  

In particular, of course, we have 
\[
\stable{\Gamma}(A,A)^{*} \cong \pi_{*}^{A}(A) \cong
\Ext^{*}_{\Gamma}(A,A).  
\]
This is the stable homotopy of the sphere in $\stable{\Gamma}$.  

The following point is also valuable. 

\begin{proposition}\label{prop-Brown}
Suppose $E$ is a commutative ring spectrum that is Landweber exact over
$BP$.  Then $\stable{E_{*}E}$ is a Brown category, so that every
homology functor is representable.  
\end{proposition}

This proposition follows from Theorem~4.1.5 of~\cite{hovey-axiomatic}.
Indeed, we can assume $E=E(n)$ or $BP$, and then one can easily check
using the cobar resolution that $\Ext^{**}_{E_{*}E}(E_{*},E_{*})$ is
countable.  

\subsection*{Ordinary homology}

We now describe ordinary homology in $\stable{\Gamma}$.  Since homotopy
isomorphisms are in particular homology isomorphisms, the ordinary
homology of a chain complex $X$ is a homology theory on
$\stable{\Gamma}$.

\begin{proposition}\label{prop-homology}
Let $\hopfalg$ be an Adams Hopf algebroid.  Ordinary homology is
represented on $\stable{\Gamma }$ by $\Gamma$ itself, as usual thought of
as a complex concentrated in degree $0$.  
\end{proposition}

Because of this proposition, we will sometimes denote $\Gamma$ by $H$.  

\begin{proof}
For a complex $X$, we have
\[
\Gamma_{*}(X) \cong \pi_{*}(Q\Gamma \smash QX),
\]
where $Q$ denotes cofibrant replacement.  Since $QX\smash (-)$ preserves
homotopy isomorphisms by Theorem~\ref{thm-model}, we have
\[
\pi_{*}(Q\Gamma \smash QX) \cong \pi_{*}(\Gamma \smash QX).
\]
Since $\Gamma \smash QX$ is already fibrant, as it is a complex of
relative injectives, we have 
\[
\pi_{*}(\Gamma \smash QX) \cong \Ch{\Gamma}(A,\Gamma \smash
QX)/\sim, 
\]
where $\sim$ denotes the chain homotopy relation.   This
is in turn isomorphic to 
\[
\Ch{A}(A, QX)/\sim\cong H_{*}(QX)\cong H_{*}X
\]
by adjointness.
\end{proof}

Note that the Hopf algebroid $(H_{*}, H_{*}H)$ associated to homology is
isomorphic to $(A, \Gamma)$ itself, concentrated in degree $0$.  Thus
$H_{*}X$ is naturally a graded $\hopfalg$-comodule, which is bigraded in
case $\hopfalg$ is graded.  We get an Adams-Novikov spectral sequence
based on $H$ whose $E_{2}$-term is
\[
E_{2}^{s,t} \cong \Ext_{\Gamma}^{s}(A, H_{t}X),
\]
which in good cases will converge to $\pi_{*}X$.  In particular, if
$X=S$, this spectral sequence is concentrated in degrees $(s,0)$, and so
collapses and converges to
\[
\pi_{*}S \cong \Ext^{*}_{\Gamma}(A,A).  
\]
Thus, if we take $\hopfalg =(BP_{*},BP_{*}BP)$, we have built a stable
homotopy category in which the usual Adams-Novikov spectral sequence
collapses.

Note that ordinary cohomology is somewhat complicated.  This is actually
already true in the derived category $\cat{D}(A)$.  Indeed, in
$\cat{D}(A)$ we have 
\[
H^{*}(X) \cong \Ch{A}(QX,S^{0}A)^{*} 
\]
and there is no really convenient interpretation of these groups.
Similarly, in $\Ch{\Gamma}$, we have
\[
H^{*}(X) \cong \Ch{\Gamma}(QX, S^{0}\Gamma)^{*} \cong \Ch{A}(QX,
S^{0}A)^{*}.  
\]

The ordinary derived category of $\Gamma$, obtained by inverting the
homology isomorphisms, is the Bousfield localization of
$\stable{\Gamma}$ with respect to $H$.  As we have said before, this is
a non-trivial localization.  Indeed, suppose $x$ is a non-nilpotent
class in $\Ext^{s}(A,A)$ with $s>0$.  Then $x$ corresponds to a self-map
$S^{-s}\xrightarrow{}S$, which is necessarily $0$ on homology.  Hence
the telescope $x^{-1}S$ will have no homology, but will be nonzero.  In
particular, if $\hopfalg =(BP_{*},BP_{*}BP)$, there are non-nilpotent
classes in $\Ext$.  For example, $\overline{\alpha_{1}}$ is
non-nilpotent when $p=2$ by~\cite[Theorem~4.4.37]{ravenel}.

\subsection*{Homology with coefficients}

We now consider ordinary homology with coefficients in an $A$-module
$B$.  

\begin{proposition}\label{prop-homology-coeffs}
Let $\hopfalg$ be an Adams Hopf algebroid, and let $B$ be an
$A$-module.  Then $\Gamma \otimes_{A}B$ represents the homology theory
$HB$ on $\stable{\Gamma}$ defined by 
\[
(HB)_{*}(X) \cong H_{*}(B\otimes_{A}QX).  
\]
If $B$ is Landweber exact over $A$, then 
\[
(HB)_{*}(X) \cong H_{*}(B\otimes_{A}X) \cong B\otimes_{A} H_{*}X.
\]
In particular, in this case the Hopf algebroid $(HB_{*}, HB_{*}HB)$ is
isomorphic to $(B, \Gamma_{B})$ concentrated in degree $0$.  
\end{proposition}

\begin{proof}
For an object $X$ of $\stable{\Gamma}$, we have 
\[
(HB)_{*}(X) \cong \pi_{*}(\Gamma \otimes_{A}B)\smash QX)\cong
H_{*}(B\otimes_{A} QX),
\]
where we have used Lemma~\ref{lem-tensor} to manipulate the tensor
product.  In particular, $HB_{*}(S)\cong B$ concentrated in degree $0$.
If $B$ is Landweber exact over $A$, then $B\otimes_{A} (-)$ will
preserve homology isomorphisms of complexes of comodules, so
\[
(HB)_{*}(X)\cong H_{*}(B\otimes_{A} X)\cong B\otimes_{A} H_{*}X.
\]
In particular, 
\[
(HB)_{*}(HB)\cong B\otimes_{A} \Gamma \otimes_{A} B \cong \Gamma_{B}.  
\]
Thus $(HB)_{*}X$ is naturally a graded comodule over $(B, \Gamma_{B})$,
when $B$ is Landweber exact over $A$.  
\end{proof}

Thus we get theories $HE(n)$ when $\hopfalg =(BP_{*},BP_{*}BP)$ and
$B=E(n)_{*}$.

The Adams-Novikov spectral sequence based on $HB$ when $B$ is Landweber
exact will then have $E_{2}$-term
\[
E_{2}^{s,t}\cong \Ext_{\Gamma_{B}}^{s}(B\otimes_{A} H_{t}X, B\otimes_{A}
H_{t}Y).  
\]
In particular, when $X=Y=S$, this spectral sequence must collapse, since
the $E_{2}$-term is concentrated where $t=0$.  However, it is not
entirely clear to what it converges.  The obvious guess is
$\pi_{*}L_{HB}S$, where $L_{HB}$ denotes Bousfield localization with
respect to $HB$.  Bousfield's convergence results~\cite{bousfield-spec}
should be re-examined to see if they apply in a more general setting to
answer this question.   

Note that if $B$ is an $A$-algebra that is also a field, then $HB$ will
be a field object of $\stable{\Gamma}$.  In particular, if $\ideal{p}$
is a prime ideal in $A$ with residue field $k_{\ideal{p}}$, then we can
form $Hk_{\ideal{p}}$.  If we apply this to the case $\hopfalg =(BP_{*},
BP_{*}BP)$, we get field spectra $HK(n)$ corresponding to the Morava
$K$-theories, but we also get many other field spectra, including $H\Fp$
corresponding to the prime ideal $I_{\infty}$.  Note that the objects
$HK(n)$ do not detect nilpotence in $\stable{BP_{*}BP}$, since there are
non-nilpotent self-maps of $S$ that are zero on homology with any
coefficients.  


\section{Landweber exactness and the stable homotopy
category}\label{sec-comparison} 

Recall that in Section~\ref{sec-Landweber} we showed that the abelian
category of $E(n)_{*}E(n)$-comodules is a localization of the abelian
category of $BP_{*}BP$-comodules.  In Section~\ref{sec-derived}, we
introduced stable homotopy categories of $E(n)_{*}E(n)$ and
$BP_{*}BP$-comodules.  It is therefore natural to conjecture that
$\stable{E(n)_{*}E(n)}$ is a Bousfield localization of
$\stable{BP_{*}BP}$.  The goal of this section is to prove this
conjecture, thereby proving Theorem~\ref{main-C}.  

\subsection*{The functor $\Phi_{*}$}\label{subsec-Phi-preserves}

Throughout this section, we let $\hopfalg =(BP_{*},BP_{*}BP)$,
$B=E(n)_{*}$, and we let $\Phi \mathcolon \hopfalg \xrightarrow{}(B,
\Gamma_{B})$ be the induced map of Hopf algebroids.

The map of Hopf algebroids $\Phi $ induces a
functor 
\[
\Phi_{*} \mathcolon \Gamma \comod \xrightarrow{} \Gamma_{B}\comod 
\]
and a left Quillen functor 
\[
\Phi_{*}\mathcolon \Ch{\Gamma }\xrightarrow{}\Ch{\Gamma_{B}}
\]
by Theorem~\ref{thm-model-natural}.  To prove Theorem~\ref{main-C}, we
must show that $\Phi_{*}$ induces a Quillen equivalence upon suitably
localizing $\Ch{\Gamma}$.  The object of the present section is to
prove the following theorem.  

\begin{theorem}\label{thm-Phi-preserves}
The functor $\Phi_{*}\mathcolon
\Ch{\Gamma}\xrightarrow{}\Ch{\Gamma_{B}}$ preserves weak equivalences.
Its right adjoint $\Phi^{*}$ reflects weak equivalences.
\end{theorem}

We prove this theorem in a series of propositions.  

\begin{proposition}\label{prop-D-closed}
Let $\cat{D}$ denote the class of all $X\in \Ch{\Gamma}$ such that the
map $\Phi_{*}QX\xrightarrow{}\Phi_{*}X$ is a weak equivalence in
$\Ch{\Gamma_{B}}$, where $Q$ is a cofibrant replacement functor in
$\Ch{\Gamma }$.  Then $\cat{D}$ is a thick subcategory.  
\end{proposition}

\begin{proof}
Note that $\cat{D}$ is obviously closed under retracts.  To see that
$\cat{D}$ is thick, suppose we have a short exact sequence 
\[
X'\xrightarrow{}X\xrightarrow{} X''
\]
in $\Ch{\Gamma}$ such that two out of three terms are in $\cat{D}$.  By
Proposition~\ref{prop-comodule-stable}(a), this is a cofiber sequence in
$\stable{\Gamma}$.  Since $\Phi_{*}Q$ is the total left derived functor
of the left Quillen functor $\Phi_{*}$, we conclude that 
\[
\Phi_{*}QX' \xrightarrow{} \Phi_{*}QX \xrightarrow{} \Phi_{*}QX''
\]
is a cofiber sequence in $\stable{\Gamma_{B}}$.  On the other hand,
because $\Phi_{*}$ is exact, the sequence 
\[
\Phi_{*}X' \xrightarrow{} \Phi_{*}X \xrightarrow{} \Phi_{*}X''
\]
is a short exact sequence in $\Ch{\Gamma_{B}}$, and hence, applying
Proposition~\ref{prop-comodule-stable}(a) again, is also a cofiber
sequence in $\stable{\Gamma_{B}}$.  There is a map from the first of
these cofiber sequences to the second, and by assumption it is an
isomorphism on two out of three terms.  Since $\stable{\Gamma_{B}}$ is a
triangulated category, we conclude that it is also an isomorphism on the
third term, and so $\cat{D}$ is thick.  
\end{proof}

Our next goal is to show that $\cat{D}$ is closed under filtered
colimits.  For this we need to recall some standard model category
theory.  Suppose $\cat{I}$ is a small category, and $\cat{M}$ is a
cofibrantly generated model category, such as $\Ch{\Gamma}$.  Then there
is a cofibrantly generated model category structure on the diagram
category $\cat{M}^{\cat{I}}$~\cite[Theorem~12.7.1]{hirschhorn} in which
the weak equivalences and fibrations are taken objectwise. Furthermore,
the cofibrations in $\cat{M}^{\cat{I}}$ are in particular objectwise
cofibrations~\cite[Proposition~12.7.3]{hirschhorn}.  

\begin{proposition}\label{prop-D-closed-filtered}
The class $\cat{D}$ of Proposition~\ref{prop-D-closed} is closed under
filtered colimits.  
\end{proposition}

\begin{proof}
Suppose $F\mathcolon \cat{I}\xrightarrow{}\Ch{\Gamma}$ is a functor from
a filtered small category $\cat{I}$ such that $F(i)\in \cat{D}$ for all
$i\in \cat{I}$.  We must show that $\colim F(i)\in \cat{D}$.  Let $QF$
be a cofibrant replacement of $F$ in the model category on
$\Ch{\Gamma}^{\cat{I}}$ discussed prior to this proposition.  Because
the constant diagram functor obviously preserves fibrations and trivial
fibrations, the colimit is a left Quillen
functor~\cite[Theorem~12.7.9]{hirschhorn}.  Hence $\colim QF$ is
cofibrant.  Furthermore, each map $QF(i)\xrightarrow{}F(i)$ is a
homotopy isomorphism, and so, since homotopy commutes with filtered
colimits, we conclude that the map $\colim QF\xrightarrow{}\colim F$ is
a weak equivalence.  Therefore, $\colim QF$ is a cofibrant replacement
of $\colim F$.  To show that $\colim F\in \cat{D}$, then, we need only
show that the map $\Phi_{*}(\colim QF)\xrightarrow{}\Phi_{*}(\colim F)$
is a homotopy isomorphism.  Since $\Phi_{*}$ itself commutes with
colimits, this is equivalent to showing that the map $\colim
\Phi_{*}QF\xrightarrow{}\colim \Phi_{*}F$ is a homotopy isomorphism.
Since $QF$ is cofibrant, and cofibrations of diagrams are in particular
objectwise cofibrations, we conclude that $QF(i)$ is a cofibrant
replacement for $F(i)$ for all $i\in \cat{I}$.  Since $F(i)\in \cat{D}$,
then, each map $\Phi_{*}QF(i)\xrightarrow{}\Phi_{*}F(i)$ is a homotopy
isomorphism.  Hence, again using the fact that homotopy commutes with
filtered colimits, $\colim \Phi_{*}QF\xrightarrow{}\colim \Phi_{*}F$ is
a homotopy isomorphism, so $\colim F\in \cat{D}$.
\end{proof}

We now know that $\cat{D}$ is a thick subcategory that is closed under
filtered colimits and (obviously) contains all the cofibrant objects of
$\Ch{\Gamma}$.  This should mean that it has to be all of
$\Ch{\Gamma}$, and that is what we now prove. 

\begin{proposition}\label{prop-Phi-Q}
If $X\in \Ch{\Gamma}$, then the map $\Phi_{*}QX\xrightarrow{}\Phi_{*}X$
is a weak equivalence.
\end{proposition}

\begin{proof}
The proposition is just saying that the class $\cat{D}$ of
Propositions~\ref{prop-D-closed} and~\ref{prop-D-closed-filtered} is all
of $\Ch{\Gamma}$.  We prove this in three steps.  We first show that the
complexes $S^{n}M$ are in $\cat{D}$, where $M$ is a finitely presented
$\Gamma$-comodule and $S^{n}M$ denotes the complex whose only non-zero
entry is $M$ in degree $n$.  We then show that all finitely presented
complexes are in $\cat{D}$, and finally, we show that every complex is a
filtered colimit of finitely presented complexes, so is in $\cat{D}$ by
Proposition~\ref{prop-D-closed-filtered}.  

For the first step, it is clear that $S^{n}A$ is in $\cat{D}$ since it
is cofibrant.  The collection of all $M$ such that $S^{n}M$ is in
$\cat{D}$ is a thick subcategory by Proposition~\ref{prop-D-closed}; by
induction, therefore, it contains $A/I_{k}$ for all $k$.  The Landweber
filtration theorem then implies that it contains all finitely presented
$M$.  

Now suppose $X$ is a finitely presented complex.  For the purposes of
the present proof, we take this to mean that $X_{n}$ is finitely
presented for all $n$ and $0$ for almost all $n$; this is in fact
equivalent to $X$ being a finitely presented object of $\Ch{\Gamma}$ in
the categorical sense.  We easily prove by induction on the nunber of
non-zero entries in $X$ that $X\in \cat{D}$.  Indeed, the base case of
one non-zero entry is handled in the preceding paragraph.  For the
induction step, let $X'$ be the subcomplex of $X$ obtained by removing
the non-zero entry in the largest possible degree.  Then $X'\in \cat{D}$
by the induction hypothesis, and the quotient $X/X'\in \cat{D}$ by the
preceding paragraph.  Since $\cat{D}$ is thick, $X\in \cat{D}$.  

Now suppose $X$ is an arbitrary complex.  Let $\cat{F}/X$ denote the
category of all maps $F\xrightarrow{}X$, where $F$ is a finitely
presented complex.  This is easily seen to have a small skeleton and to
be a filtered category.  There is an obvious inclusion functor
$i\mathcolon \cat{F}/X\xrightarrow{}\Ch{\Gamma}$, and an obvious map
\[
f\mathcolon \colim i\xrightarrow{}X.
\]
We claim that $f$ is an isomorphism.  To see that $f$ is surjective,
choose $x\in X_{n}$.  Since the comodule $X_{n}$ is a filtered colimit
of finitely presented comodules, there is a finitely presented comodule
$F$ and a map $F\xrightarrow{}X_{n}$ whose image contains $x$.  This
gives a map of complexes $D^{n}F\xrightarrow{}X$ whose image contains
$x$, and so $f$ is surjective.  To see that $f$ is injective, suppose
$j\mathcolon F\xrightarrow{}X$ is an object of $\cat{F}/X$ and $x\in
F_{n}$ has $jx=0$.  Let $K$ denote the kernel of the map $j$, so that
$x\in K_{n}$.  Now $K$ may not be finitely presented, but at least there
is a map $F'\xrightarrow{}K_{n}$ from a finitely presented comodule
whose image contains $x$.  This corresponds to a map
$D^{n}F'\xrightarrow{}K$ of complexes, which induces an object
$F/D^{n}F'\xrightarrow{}X$ of $\cat{F}/X$ and a map
$F\xrightarrow{}F/D^{n}F'$in $\cat{F}/X$ that sends $x$ to $0$.  Thus
$x$ is $0$ in $\colim i$ and so $f$ is injective.  
\end{proof}

We can now give the proof of Theorem~\ref{thm-Phi-preserves}.

\begin{proof}[Proof of Theorem~\ref{thm-Phi-preserves}]
Suppose $f\mathcolon X\xrightarrow{}Y$ is a weak equivalence in
$\Ch{\Gamma}$.  Then $Qf$ is a weak equivalence between cofibrant
objects, so $\Phi_{*}Qf$ is a weak equivalence since $\Phi_{*}$ is a
left Quillen functor.  But we have a commutative square 
\[
\begin{CD}
\Phi_{*}QX @>\Phi_{*}Qf>> \Phi_{*}QY \\
@VVV @VVV \\
\Phi_{*}X @>>\Phi_{*}f> \Phi_{*}Y
\end{CD}
\]
where the vertical maps are weak equivalences, by
Proposition~\ref{prop-Phi-Q}.  Hence $\Phi_{*}f$ is a weak equivalence
as well.  

To prove the second part of Theorem~\ref{thm-Phi-preserves}, suppose $f$
is a map in $\Ch{\Gamma_{B}}$ such that $\Phi^{*}f$ is a weak
equivalence in $\Ch{\Gamma}$.  Then $\Phi_{*}\Phi^{*}f$ is a weak
equivalence in $\Ch{\Gamma}$ by what we have just proved.  But $f$ is
naturally isomorphic to $\Phi_{*}\Phi^{*}f$ by
Theorem~\ref{thm-Giraud}.  
\end{proof}

\subsection*{Localization}\label{subsec-comp}

We have just seen that $\Phi_{*}\mathcolon
\Ch{\Gamma}\xrightarrow{}\Ch{\Gamma_{B}}$ preserves weak equivalences.
But of course it does not reflect weak equivalences, since
$\Phi_{*}(A/I_{n+1})=0$.  We dealt with this problem already in the
abelian category world by localizing $\Gamma \comod$ so as to force
$0\xrightarrow{}A/I_{n+1}$ to be an isomorphism.  We now want to do the
same thing for chain complexes.  

More precisely, we define $L_{n}^{f}\Ch{\Gamma }$ to be the category
$\Ch{\Gamma }$ equipped with the model structure that is the Bousfield
localization of the homotopy model structure with respect to the maps
$0\xrightarrow{}s^{k}A/I_{n+1}$ for all $k$, where $s^{k}A/I_{n+1}$
denotes the complex which is $A/I_{n+1}$ in degree $k$ and $0$
elsewhere.  Recall from~\cite{hirschhorn} that this means that a left
Quillen functor
\[
F\mathcolon \Ch{\Gamma } \xrightarrow{} \cat{M}
\]
defines a left Quillen functor 
\[
F \mathcolon L_{n}^{f}\Ch{\Gamma } \xrightarrow{} \cat{M}
\]
if and only if $0\xrightarrow{}(LF)(s^{k}A/I_{n+1})$ is an isomorphism
in $\sho \cat{M}$, where $LF$ denotes the total left derived functor of
$F$.  The homotopy category of $L_{n}^{f}\Ch{\Gamma }$ is the finite
localization $L_{n}^{f}\stable{\Gamma }$ in the sense of
Miller~\cite{miller-finite} of $\stable{\Gamma }$ away from $A/I_{n+1}$.
The total left derived functor of the identity, thought of as a functor
from $\Ch{\Gamma }$ to $L_{n}^{f}\Ch{\Gamma }$, is the finite
localization functor $L_{n}^{f}$ on $\stable{\Gamma }$.

\begin{proposition}\label{prop-induced-local}
The Quillen functor 
\[
\Phi_{*}\mathcolon \Ch{\Gamma}\xrightarrow{}\Ch{\Gamma_{B}}
\]
induces a left Quillen functor
\[
\Phi_{*}\mathcolon L_{n}^{f}\Ch{\Gamma}\xrightarrow{} \Ch{\Gamma_{B}}.
\]
\end{proposition}

\begin{proof}
We need to show that $(L\Phi_{*})(A/I_{n+1})=0$.  In light of
Proposition~\ref{prop-Phi-Q}, $(L\Phi_{*})X$ is naturally isomorphic to
$\Phi_{*}X$ for any $X\in \Ch{\Gamma}$.  Thus
\[
(L\Phi_{*})(A/I_{n+1})\cong \Phi_{*}(A/I_{n+1})=0. \ \ \Box 
\]
\renewcommand{\qed}{}\end{proof}

Note that $\Phi_{*}$ still preserves weak equivalences when thought as a
functor from $L_{n}^{f}\Ch{\Gamma}$,

\begin{proposition}\label{prop-induced-preserves}
The Quillen functor 
\[
\Phi_{*}\mathcolon L_{n}^{f}\Ch{\Gamma}\xrightarrow{}\Ch{\Gamma_{B}}
\]
preserves weak equivalences, and its right adjoint $\Phi^{*}$
reflects weak equivalences.  
\end{proposition}

\begin{proof}
Suppose $f\mathcolon X\xrightarrow{}Y$ is a weak equivalence.  Factor
$f=pi$, where $i$ is a trivial cofibration and $p$ is a trivial
fibration.  Then $\Phi_{*}i$ is a weak equivalence since $\Phi_{*}$ is a
left Quillen functor.  On the other hand, since the trivial fibrations
do not change under Bousfield localization, $p$ is a weak equivalence in
$\Ch{\Gamma}$.  Thus Theorem~\ref{thm-Phi-preserves} implies that
$\Phi_{*}p$ is a weak equivalence.  Hence
$\Phi_{*}f=(\Phi_{*}p)(\Phi_{*}i)$ is a weak equivalence, as required.  

The proof that $\Phi^{*}$ reflects weak equivalences is the same as the
proof of the corresponding part of Theorem~\ref{thm-Phi-preserves}.  
\end{proof}

To prove Theorem~\ref{main-C}, we will show that 
\[
\Phi_{*}\mathcolon L_{n}^{f}\Ch{\Gamma} \xrightarrow{} \Ch{\Gamma_{B}}
\]
is a Quillen equivalence.  To do so, we will use the following lemma,
which is proved in~\cite[Corollary~1.3.16]{hovey-model}.

\begin{lemma}\label{lem-Quillen}
Suppose $F\mathcolon \cat{C}\xrightarrow{}\cat{D}$ is a left Quillen
functor of model categories, with right adjoint $U$.  Then $F$ is a
Quillen equivalence if and only if the following two conditions hold.
\begin{enumerate}
\item [(a)] $U$ reflects weak equivalences between fibrant objects. 
\item [(b)] The map $X\xrightarrow{}URFX$ is a weak equivalence for all
cofibrant $X$ in $\cat{C}$, where $R$ denotes fibrant replacement and
the map is induced by the unit of the adjunction.  
\end{enumerate}
\end{lemma}

We have already seen in Proposition~\ref{prop-induced-preserves} that
$\Phi^{*}$ reflects all weak equivalences.  We point out that there is a
much simpler proof that $\Phi_{*}$ reflects weak equivalences between
fibrant objects; if $X$ is fibrant in $\Ch{\Gamma_{B}}$, then
adjointness implies that $\pi_{*}(\Phi^{*}X)\cong \pi_{*}X$.  

The other condition of Lemma~\ref{lem-Quillen} is harder to check.  Here
are the main points of the argument.  
\begin{itemize}\label{argument}
\item [(a)] We first show, using the fact that $\Phi^{*}$ preserves
filtered colimits, that it suffices to show that
$A\xrightarrow{}\Phi^{*}(LB)$ is an $L_{n}^{f}$-equivalence, where $LB$
denotes the cobar resolution of $B$.  
\item [(b)] A Bousfield class argument that shows that it suffices to
prove that 
\[
v_{k}^{-1}A/I_{k}\xrightarrow{}v_{k}^{-1}A/I_{k}\smash Q\Phi^{*}(LB)
\]
is a homotopy isomorphism for all $0\leq k\leq n$, where $Q$ denotes
cofibrant replacement.
\item [(c)] We show that, although $\Phi^{*}(LB)$ is not cofibrant, it
is still nice enough that it suffices to check that
\[
v_{k}^{-1}A/I_{k} \xrightarrow{} v_{k}^{-1}A/I_{k} \smash \Phi^{*}(LB)
\]
is a homotopy isomorphism.  
\item [(d)] It was proved in~\cite{hovey-hopf} that the Hopf algebroid
$(v_{k}^{-1}A/I_{k},v_{k}^{-1}\Gamma/I_{k})$ is weakly equivalent to
$(v_{k}^{-1}B/I_{k}, v_{k}^{-1}\Gamma_{B}/I_{k})$.  Hence it suffices to
show that 
\[
v_{k}^{-1}B/I_{k} \xrightarrow{} v_{k}^{-1}B/I_{k} \smash
\Phi_{*}\Phi^{*}(LB) \cong v_{k}^{-1}B/I_{k}\smash LB
\]
is a homotopy isomorphism, and this is clear.  
\end{itemize}

We now fill in the details of this argument, beginning with Step~(a).

\begin{proposition}\label{prop-reduction}
The Quillen functor 
\[
\Phi_{*}\mathcolon L_{n}^{f}\Ch{\Gamma} \xrightarrow{} \Ch{\Gamma_{B}}
\]
is a Quillen equivalence if and only if the map 
\[
A \xrightarrow{} \Phi^{*}(LB)
\]
is an $L_{n}^{f}$-equivalence.  
\end{proposition}

Recall that $LB$ denotes the cobar resolution of $B$ as a
$\Gamma_{B}$-comodule.  

Before proving this proposition, we need a lemma.

\begin{lemma}\label{lem-right-coproduct}
The total right derived functor 
\[
R\Phi^{*} \mathcolon \stable{\Gamma_{B}} \xrightarrow{}
L_{n}^{f}\stable{\Gamma}
\]
preserves coproducts.  
\end{lemma}

\begin{proof}
First note that because $L_{n}^{f}$ is a finite localization, it is in
particular smashing~\cite{miller-finite}.  Thus the coproduct in
$L_{n}^{f}\stable{\Gamma}$ is the same as the coproduct in
$\stable{\Gamma}$.  Also note that a fibrant replacement in
$\stable{\Gamma_{B}}$ is given by $LB\smash (-)$, which clearly
preserves coproducts.  Hence, it suffices to show that 
\[
\Phi^{*} \mathcolon \Gamma_{B}\comod \xrightarrow{} \Gamma \comod 
\]
preserves coproducts.  Since $\Phi^{*}$ certainly preserves finite
coproducts, and any coproduct is a filtered colimit of finite
coproducts, it suffices to show that $\Phi^{*}$ preserves all filtered
colimits.  This follows from the fact that $L_{n}=\Phi^{*}\Phi_{*}$
preserves filtered colimits (Theorem~\ref{thm-localization}); more
details can be found in~\cite{hovey-strickland-derived}.
\end{proof}

\begin{proof}[Proof of Proposition~\ref{prop-reduction}]
Lemma~\ref{lem-Quillen} tells us that if $\Phi_{*}$ is a Quillen
equivalence, then $A\xrightarrow{}\Phi^{*}R\Phi_{*}A$ must be a weak
equivalence in $L_{n}^{f}\Ch{\Gamma}$.  Since $\Phi_{*}A=B$, and since
$LB$ is a fibrant replacement for $B$ in $\Ch{\Gamma_{B}}$, we conclude
that $A\xrightarrow{}\Phi^{*}(LB)$ is an $L_{n}^{f}$-equivalence.  

Conversely, suppose $A\xrightarrow{}\Phi^{*}(LB)$ is an
$L_{n}^{f}$-equivalence.  By Lemma~\ref{lem-Quillen} and
Proposition~\ref{prop-induced-preserves}, it suffices to show that
$X\xrightarrow{}\Phi^{*}R\Phi_{*}X$ is an $L_{n}^{f}$-equivalence for
all cofibrant $X$.  This is equivalent to proving that 
\[
X \xrightarrow{\eta_{X} } (R\Phi^{*})(L\Phi_{*})X
\] 
is an isomorphism in $L_{n}^{f}\stable{\Gamma}$ for all $X$, where
$R\Phi^{*}$ denotes the total right derived functor of $\Phi^{*}$ and
$L\Phi_{*}$ denotes the total left derived functor of $\Phi_{*}$.  Let
$\cat{D}$ denote the full subcategory of $L_{n}^{f}\stable{\Gamma }$ of
those $X$ such that $\eta_{X}$ is an isomorphism.  By hypothesis,
$\cat{D}$ contains $L_{n}^{f}A$.  Since $L\Phi_{*}$ and $R\Phi^{*}$ both
preserve exact triangles, $\cat{D}$ is a thick subcategory.  As a left
adjoint, $L\Phi_{*}$ preserves coproducts, and
Lemma~\ref{lem-right-coproduct} assures us that $R\Phi^{*}$ also
preserves coproducts.  Thus $\cat{D}$ is a localizing subcategory.  In
any monogenic stable homotopy category, the only localizing subcategory
that contains the unit is the whole category.
\end{proof}

We are now reduced to showing that $A\xrightarrow{}\Phi^{*}(LB)$ is an
$L_{n}^{f}$-equivalence.  The theory of Bousfield classes gives us the
following lemma.  

\begin{lemma}\label{lem-Bousfield}
A map $f$ of cofibrant objects in $\Ch{\Gamma}$ is an
$L_{n}^{f}$-equivalence if and only if $v_{k}^{-1}A/I_{k}\smash f$ is a
weak equivalence for all $k$ with $0\leq k\leq n$.  
\end{lemma}

\begin{proof}
As usual, let $\Bousfield{X}$ denote the Bousfield class of $X$ in
$\stable{\Gamma }$, thought of as the collection of all $Y$ in
$\stable{\Gamma}$ such that $X\smash Y=0$.  The cofiber sequences
\[
A/I_{k} \xrightarrow{v_{k}} A/I_{k} \xrightarrow{} A/I_{k+1}
\]
imply that 
\[
\Bousfield{A/I_{k}} = \Bousfield{v_{k}^{-1}A/I_{k}} \vee
\Bousfield{A/I_{k+1}},
\]
by~\cite[Lemma~1.34]{rav-loc}.  Thus, we find 
\[
\Bousfield{A} = \bigvee_{k=0}^{n} \Bousfield{v_{k}^{-1}A/I_{k}} \vee
\Bousfield{A/I_{n+1}}.  
\]
As in the usual stable homotopy category, this implies that $L_{n}^{f}$
is localization with respect to $\bigoplus_{k=0}^{n}
v_{k}^{-1}A/I_{k}$.  

It follows that $f$ is an $L_{n}^{f}$-equivalence if and only if
$v_{k}^{-1}A/I_{k} \smash^{L} f$ is a weak equivalence for all $k$ with
$0\leq k\leq n$, where $\smash^{L}$ denotes the total left derived
functor of $\smash$.  Recall from Theorem~\ref{thm-model} that if $X$ is
cofibrant, then $X\smash (-)$ preserves homotopy isomorphisms.  It
follows that, if $X$ is cofibrant, $(-)\smash^{L}X\cong (-)\smash X$ in
$\stable{\Gamma}$.  It follows that, if $f$ is a map of cofibrant
objects, then $f$ is an $L_{n}^{f}$-equivalence if and only if
$v_{k}^{-1}A/I_{k}\smash f$ is a weak equivalence for all $k$ with
$0\leq k\leq n$.
\end{proof}

By combining Proposition~\ref{prop-reduction} with
Lemma~\ref{lem-Bousfield}, we get the following corollary, which is
Step~(b) of the argument on page~\pageref{argument}.  

\begin{corollary}\label{cor-Bousfield}
The Quillen functor 
\[
\Phi_{*}\mathcolon L_{n}^{f}\Ch{\Gamma} \xrightarrow{} \Ch{\Gamma_{B}}
\]
is a Quillen equivalence if and only if the map 
\[
v_{k}^{-1}A/I_{k} \xrightarrow{} v_{k}^{-1}A/I_{k} \smash Q\Phi^{*}(LB)
\]
is a homotopy isomorphism for all $0\leq k\leq n$, where $Q$ denotes a
cofibrant replacement functor in $\Ch{\Gamma}$.  
\end{corollary}

To accomplish Step~(c) of the argument on page~\pageref{argument}, we
need to know something about $\Phi^{*}(LB)$.

\begin{lemma}\label{lem-no-tor}
We have 
\[
\Tor_{A}^{j}(A/I_{k}, \Phi^{*}(LB)_{m})=0
\]
for all $j>0$, $k\geq 0$, and $m\in \Z$.  
\end{lemma}

\begin{proof}
Recall that 
\[
(LB)_{-m}=\Gamma_{B}\otimes_{B} \overline{\Gamma_{B}}^{\otimes_{B}m}
\]
for $m\geq 0$, and is $0$
otherwise.  Here $\overline{\Gamma_{B}}$ is the cokernel of the left
unit $\eta_{L}\mathcolon B\xrightarrow{}\Gamma_{B}$.  Since this map is
split as a map of $B$-modules, $\overline{\Gamma_{B}}$ is a flat
$B$-module.  Therefore $\overline{\Gamma_{B}}^{\otimes_{B}m}$ is also a
flat $B$-module, and hence a filtered colimit of projective $B$-modules.  

Now, we have 
\[
\Phi^{*}(LB)_{-m} = \Gamma \otimes_{A}
\overline{\Gamma_{B}}^{\otimes_{B}m}.
\]
Since $\Tor^{j}_{A}(A/I_{k},-)$ commutes with filtered colimits, it suffices
to show that 
\[
\Tor_{A}^{j}(A/I_{k},\Gamma \otimes_{A} M)=0
\]
for all $j>0$, all $k$, and all projective $B$-modules $M$.  But then we
can easily reduce to the case $M=B$, so we must show that 
\[
\Tor_{A}^{j}(A/I_{k}, \Gamma \otimes_{A}B) =0
\]
for all $j>0$ and all $k$.  We prove this by induction on $k$, using the
exact sequences 
\[
0 \xrightarrow{} A/I_{k} \xrightarrow{v_{k}} A/I_{k} \xrightarrow{}
A/I_{k+1} \xrightarrow{} 0.  
\]
We are reduced to showing that $v_{k}$ is not a zero-divisor on $(\Gamma
\otimes_{A}B)/I_{k}$.  Since $I_{k}$ is invariant, this is the same as
showing that $v_{k}$ is not a zero-divisor on $\Gamma
\otimes_{A}(B/I_{k})$.  Since $v_{k}$ is itself invariant modulo $I_{k}$
and $\Gamma$ is flat over $A$, this is in turn equivalent to showing
that $v_{k}$ is not a zero-divisor on $B/I_{k}$.  This follows because
$B$ is Landweber exact.  
\end{proof}

With this lemma in hand, we can carry out Step~(c) of our argument.  

\begin{proposition}\label{prop-dont-need-Q}
The map 
\[
v_{k}^{-1}A/I_{k}\smash
Q\Phi^{*}(LB)\xrightarrow{}v_{k}^{-1}A/I_{k} \smash \Phi^{*}(LB)
\]
is a homotopy isomorphism in $\Ch{\Gamma }$.  In particular, 
\[
\Phi_{*} \mathcolon L_{n}^{f}\Ch{\Gamma} \xrightarrow{} \Ch{\Gamma_{B}}
\]
is a Quillen equivalence if and only if 
\[
v_{k}^{-1}A/I_{k} \xrightarrow{} v_{k}^{-1}A/I_{k} \smash \Phi^{*}(LB)
\]
is a homotopy isomorphism for all $0\leq k\leq n$.  
\end{proposition}

\begin{proof}
Let $q$ denote the map
\[
q\mathcolon Q\Phi^{*}(LB)\xrightarrow{}\Phi^{*}(LB).
\]
Because homotopy commutes with filtered colimits, it suffices to show
that $A/I_{k}\smash q$ is a homotopy isomorphism for all $k$.  

Now, $q$ is a trivial fibration in the homotopy model structure, so we
have a short exact sequence of complexes
\[
0 \xrightarrow{} K \xrightarrow{} Q\Phi^{*}(LB) \xrightarrow{q}
\Phi^{*}(LB) \xrightarrow{} 0.
\]
The long exact sequence in homotopy implies that $\pi_{t}^{A}(K)=0$ for
all $n$.  Lemma~\ref{lem-no-tor} implies that we have a short exact
sequence 
\[
0 \xrightarrow{} A/I_{k}\smash K \xrightarrow{}A/I_{k} \smash
Q\Phi^{*}(LB) \xrightarrow{} A/I_{k} \smash \Phi^{*}(LB) \xrightarrow{} 0
\]
for all $k$.  The long exact sequence in homotopy implies that we need
only check that $\pi_{t}^{A}(A/I_{k}\smash K)=0$ for all $n$.  

To se this, note that the short exact sequence defining $K$ realizes
$K_{m}$ as the first syzygy of $\Phi^{*}(LB)_{m}$, since
$(Q\Phi^{*}(LB))_{m}$ is projective over $A$.  Hence
Lemma~\ref{lem-no-tor} implies that
\[
\Tor_{A}^{j}(A/I_{k},K_{m})=0
\]
for all $j>0$.  Hence we have short exact sequences
\[
0 \xrightarrow{} A/I_{k} \smash K \xrightarrow{} A/I_{k}\smash K
\xrightarrow{} A/I_{k+1} \smash K \xrightarrow{} 0
\]
for all $k$.  The long exact sequence in homotopy and induction on $k$
now complete the proof.
\end{proof}

The final step of the argument on page~\pageref{argument} requires us to
know more about the Hopf algebroids $(v_{k}^{-1}A/I_{k},
v_{k}^{-1}\Gamma /I_{k})$ and $(v_{k}^{-1}B/I_{k},
v_{k}^{-1}\Gamma_{B}/I_{k})$.

\begin{lemma}\label{lem-localized-Hopf}
Let $(C, \Sigma)$ denote $(v_{k}^{-1}A/I_{k}, v_{k}^{-1}\Gamma /I_{k})$,
and let $(C_{B}, \Sigma_{B})$ denote $(v_{k}^{-1}B/I_{k},
v_{k}^{-1}\Gamma /I_{k})$.  Then\uc 
\begin{enumerate}
\item [(a)] Both $(C, \Sigma)$ and $(C_{B}, \Sigma_{B})$ are Adams Hopf
algebroids\usc  
\item [(b)] The stable homotopy categories $\stable{\Sigma}$ and
$\stable{\Sigma_{B}}$ are monogenic\usc 
\item [(c)] If $0\leq k\leq n$, $\Phi$ induces a weak equivalence of
Hopf algebroids 
\[
\Phi \mathcolon (C,\Sigma)\xrightarrow{}(C_{B}, \Sigma_{B}).
\]
\end{enumerate}
\end{lemma}

\begin{proof}
For part~(a), we know that $\Gamma$ is a filtered colimit of comodules
that are finitely generated projective $A$-modules.  By tensoring with
$v_{k}^{-1}A/I_{k}$, we se that $\Sigma$ is a filtered colimit of
comodules that are finitely generated projective $C$-modules, and so
$(C, \Sigma)$ is an Adams Hopf algebroid.  Similarly, $(C_{B},
\Sigma_{B})$ is an Adams Hopf algebroid.  

For part~(b), the proof is again the same for $(C, \Sigma)$ and $(C_{B},
\Sigma_{B})$, so we concentrate on $(C, \Sigma)$.  We will use
Proposition~\ref{prop-monogenic}, so we need to show that every
dualizable $\Sigma$-comodule is in the thick subcategory generated by
$C$.  We will do this by showing that every finitely presented
$\Sigma$-comodule has a Landweber filtration.  To do so, we will use the
Hopf algebroid $(v_{k}^{-1}A, v_{k}^{-1}\Gamma v_{k}^{-1})$, obtained by
inverting $v_{k}$ and $\eta_{R}v_{k}$.  A $\Sigma$-comodule $M$ is just
a $(v_{k}^{-1}\Gamma v_{k}^{-1})$-comodule on which $I_{k}$ acts
trivially.  Since $I_{k}$ is finitely generated, $M$ is
finitely presented if and only if it is finitely presented as a
$v_{k}^{-1}\Gamma v_{k}^{-1}$-comodule.  Since $v_{k}^{-1}A$ is
Landweber exact of height $k$, Theorem~\ref{thm-Landweber} gives us a
Landweber filtration of $M$ as a $v_{k}^{-1}\Gamma v_{k}^{-1}$-comodule
in which each filtration quotient is isomorphic to $v_{k}^{-1}A/I_{j}$
for some $j\leq k$.  Since $M$ is killed by $I_{k}$, in fact each
filtration quotient must be isomorphic to $v_{k}^{-1}A/I_{k}$, giving us
our Landweber filtration of $M$ as a $\Sigma$-comodule.  

Part~(c) is a special case of Theorem~E of~\cite{hovey-hopf}.  
\end{proof}

Lemma~\ref{lem-localized-Hopf} allows us to carry out the final step of
our argument.  

\begin{proposition}\label{prop-tel-equiv}
Suppose $f$ is a map in $\Ch{\Gamma}$, and $0\leq k\leq n$.  Then
$v_{k}^{-1}A/I_{k}\smash f$ is a homotopy isomorphism in $\Ch{\Gamma}$
if and only if $v_{k}^{-1}B/I_{k}\smash \Phi_{*}f$ is a homotopy
isomorphism in $\Ch{\Gamma_{B}}$.
\end{proposition}

\begin{proof}
Let $C=v_{k}^{-1}A/I_{k}$ and let $\Sigma =v_{k}^{-1}\Gamma /I_{k}$, as
in Lemma~\ref{lem-localized-Hopf}.  The category of
$\Sigma$-comodules is just the full subcategory of $\Gamma$-comodules on
which $I_{k}$ acts trivially and $v_{k}$ acts invertibly.  This follows
from the fact that $I_{k}$ is invariant and $v_{k}$ is primitive modulo
$I_{k}$.  Thus, if $X$ is a complex in $\Ch{\Sigma}$, we can also think
of $X$ as a complex in $\Ch{\Gamma}$.  As such, $X$ has homotopy groups
$\pi_{*}^{C}(X)$ and $\pi_{*}^{A}(X)$.  We claim that these are
naturally isomorphic.  Indeed, one can easily check that $LC$, the cobar
complex of $(C, \Sigma)$, is just $v_{k}^{-1}LA/I_{k}$.  Hence $LC\smash
X$ is naturally isomorphic to $LA\smash X$.  From this, one can easily
check the desired isomorphism.  

Therefore, using parts~(a) and~(b) of Lemma~\ref{lem-localized-Hopf}, we
conclude that $v_{k}^{-1}A/I_{k}\smash f$ is a homotopy isomorphism in
$\Ch{\Gamma}$ if and only if it is a weak equivalence in $\Ch{\Sigma}$.  
Similarly, let  $C_{B}=v_{k}^{-1}B/I_{k}$ and
$\Sigma_{B}=v_{k}^{-1}\Gamma_{B}/I_{k}$.  We find that
$v_{k}^{-1}B/I_{k}\smash \Phi_{*}f$ is a homotopy isomorphism in 
$\Ch{\Gamma_{B}}$ if and only if it is a weak equivalence in
$\Ch{\Sigma_{B}}$.  

Now, use part~(c) of Lemma~\ref{lem-localized-Hopf} and
Theorem~\ref{thm-model-natural} to conclude that
$v_{k}^{-1}A/I_{k}\smash f$ is a weak equivalence in $\Ch{\Sigma}$ if
and only if 
\[
v_{k}^{-1}B/I_{k} \otimes_{v_{k}^{-1}A/I_{k}} (v_{k}^{-1}A/I_{k} \smash f)
\cong  v_{k}^{-1}B/I_{k} \smash \Phi_{*}f
\]
is a weak equivalence in $\Ch{\Sigma_{B}}$.  
\end{proof}

We can now complete the proof of Theorem~\ref{main-C}, which we first
restate in stronger form.  

\begin{theorem}\label{thm-Quillen-equiv}
The Quillen functor 
\[
\Phi_{*} \mathcolon L_{n}^{f}\Ch{\Gamma} \xrightarrow{} \Ch{\Gamma_{B}}
\]
is a Quillen equivalence.  Furthermore, $\Phi_{*}$ and its right adjoint
$\Phi^{*}$ preserve and reflect weak equivalences.  
\end{theorem}

\begin{proof}
Proposition~\ref{prop-dont-need-Q} and Proposition~\ref{prop-tel-equiv}
imply that, to check that $\Phi_{*}$ is a Quillen equivalence, we only
need to check that the map
\[
v_{k}^{-1}B/I_{k} \xrightarrow{} v_{k}^{-1}B/I_{k} \smash
\Phi_{*}\Phi^{*}(LB)
\]
is a homotopy isomorphism for $0\leq k\leq n$.  But $\Phi_{*}\Phi^{*}$
is natually isomorphic to the identity functor by
Theorem~\ref{thm-Giraud}.  Hence we need only check that the map 
\[
v_{k}^{-1}B/I_{k} \xrightarrow{} v_{k}^{-1}B/I_{k} \smash LB
\]
is a homotopy isomorphism, which follows from
Proposition~\ref{prop-homotopy-groups}(d).  

We have already seen that $\Phi_{*}$ preserves weak equivalences and
that $\Phi^{*}$ reflects them in
Proposition~\ref{prop-induced-preserves}.  Suppose $f$ is a map in
$\Ch{\Gamma}$ such that $\Phi_{*}f$ is a weak equivalence.  Since
$\Phi_{*}$ preserves weak equivalences, it follows that $\Phi_{*}Qf$ is
a weak equivalence.  But, since $\Phi_{*}$ is a Quillen
equivalence, it must reflect weak equivalences between cofibrant objects
by~\cite[Corollary~1.3.16]{hovey-model}.  Hence $Qf$ is a weak
equivalence, and so $f$ is a weak equivalence.  This proves that
$\Phi_{*}$ reflects weak equivalences.  

Now suppose $g$ is a weak equivalence in $\Ch{\Gamma_{B}}$.  By
Theorem~\ref{thm-Giraud}, $g$ is naturally isomorphic to
$\Phi_{*}\Phi^{*}g$.  Since $\Phi_{*}$ reflects weak equivalences, we
conclude that $\Phi^{*}g$ is a weak equivalence.  
\end{proof}

We point out that it is possible to further localize the Quillen
equivalence in Theorem~\ref{thm-Quillen-equiv} to obtain a Quillen
equivalence 
\[
\Phi_{*}\mathcolon L_{n}\Ch{\Gamma} \xrightarrow{} L_{n}\Ch{\Gamma_{B}}
\]
where $L_{n}$ is Bousfield localization with respect to $HE(n)$ is the
first case, and ordinary homology $H$ in the second.  

\subsection*{The change of rings theorem}

In this final part, we show how our work implies the Miller-Ravenel change
of rings theorem.  To begin with, here is our generic change of rings
theorem.  Recall our notational conventions: $\hopfalg
=(BP_{*},BP_{*}BP)$, $B=E(n)_{*}$, and $\Gamma_{B}=E(n)_{*}E(n)$.  

\begin{theorem}\label{thm-change-of-ring}
Suppose $X\in \Ch{\Gamma}$ and $Y$ is an $L_{n}^{f}$-local object of
$\stable{\Gamma}$.  Then 
\[
\stable{\Gamma}(X,Y)^{*} \cong \stable{\Gamma_{B}}(\Phi_{*}X,\Phi_{*}Y)^{*}.
\] 
\end{theorem}

\begin{proof}
First of all, since $\Phi_{*}$ preserves weak equivalences, we have 
\[
\stable{\Gamma_{B}}(\Phi_{*}X, \Phi_{*}Y) \cong
\stable{\Gamma_{B}}(\Phi_{*}QX, \Phi_{*}QY).  
\]
Also, by Theorem~\ref{thm-Quillen-equiv}, we have 
\[
\stable{\Gamma_{B}}(\Phi_{*}QX, \Phi_{*}QY) \cong
(L_{n}^{f}\stable{\Gamma})(X,Y). 
\]
Since $Y$ is already $L_{n}^{f}$-local, we have 
\[
(L_{n}^{f}\stable{\Gamma})(X,Y) \cong \stable{\Gamma}(X,Y),
\]
as required.  
\end{proof}

We claim that this corollary captures the Miller-Ravenel change of rings
theorem.  To see this, we need the following lemma.

\begin{lemma}\label{lem-module-local}
Suppose $N$ is a $\Gamma$-comodule with $L_{n}N=N$ and $L_{n}^{i}N=0$
for $i>0$.  Then $N$ is an $L_{n}^{f}$-local object of
$\stable{\Gamma}$.
\end{lemma}

Given this lemma, we immediately get the following corollary, which is
Theorem~\ref{main-D} and more like the usual change of rings theorems.

\begin{corollary}\label{cor-miller-ravenel}
Suppose $M$ is a finitely presented $BP_{*}BP$-comodule, and $N$ is a
$BP_{*}BP$-comodule such that $L_{n}N=N$ and $L_{n}^{i}N=0$ for all
$i>0$.  Then 
\[
\Ext_{BP_{*}BP}^{**}(M,N) \cong
\Ext_{E(n)_{*}E(n)}(E(n)_{*}\otimes_{BP_{*}}M,
E(n)_{*}\otimes_{BP_{*}}N).  
\]
\end{corollary}

This corollary includes both the Miller-Ravenel change of rings
theorem~\cite{miller-ravenel}, by taking $N$ with $N=v_{n}^{-1}N$, and
the change of rings theorem of the author and
Sadofsky~\cite{hovey-sadofsky-picard}, by taking $N$ with
$v_{j}^{-1}N=N$ for some $j\leq n$.  Here we are using
Theorem~\ref{thm-localization}(e) to verify the hypothesis of 
Corollary~\ref{cor-miller-ravenel}.  

\begin{proof}
Simply apply Theorem~\ref{thm-change-of-ring}, using
Lemma~\ref{lem-module-local} to see that $N$ is $L_{n}^{f}$-local, and
Proposition~\ref{prop-comodule-stable}(c) to identify the groups in
question as $\Ext$ groups.  
\end{proof}

We still owe the reader a proof of Lemma~\ref{lem-module-local}, which,
incidentally, is presumably a special case of a spectral sequence that
computes $H_{*}(L_{n}^{f}X)$ for $X\in \stable{\Gamma}$ from the derived
functors $L_{n}^{i}H_{*}X$, in analogy to the spectral sequence of
Theorem~\ref{thm-spectral}.  The converse of
Lemma~\ref{lem-module-local} is true as well, though we do not need it.

\begin{proof}[Proof of Lemma~\ref{lem-module-local}]
Note that, by definition, $N$ is $L_{n}^{f}$-local if and only if
\[
\stable{\Gamma}(A/I_{n+1},N)^{*}=0.
\]
This is equivalent to 
\[
\Ext_{\Gamma}^{i}(A/I_{n+1},N)=0
\]
for all $i$, by Proposition~\ref{prop-comodule-stable}(c).  Now
$\Ext_{\Gamma}^{i}(A/I_{n+1},N)=0$ for $i=0,1$ if and only if
$L_{n}N=N$, by Theorem~\ref{thm-localization}(a).  Suppose in addition
that $L_{n}^{i}N=0$ for $i>0$.  We claim that we can find
cosyzygies $C_{j}$ of $N$ in the category of $\Gamma$-comodules such
that $L_{n}C_{j}=C_{j}$ and $L_{n}^{i}C_{j}=0$ for $i>0$.  If we
assume this, then for $i>0$ we have
\[
\Ext^{i}_{\Gamma}(A/I_{n+1},N) \cong \Ext^{1}_{\Gamma }(A/I_{n+1},
C_{i-1}) =0.
\]

We now construct the cosyzygies $C_{j}$ by induction on $j$, taking
$C_{0}=N$.  Suppose we have constructed $C_{j}$.  Since
$L_{n}C_{j}=C_{j}$, $C_{j}$ has no $v_{n}$-torsion.  It follows that the
injective hull, or indeed any essential extension of $C_{j}$, has no
$v_{n}$-torsion.  Therefore, we can find an short exact sequence 
\[
0 \xrightarrow{} C_{j} \xrightarrow{} I_{j} \xrightarrow{} C_{j+1}
\xrightarrow{} 0
\]
of $\Gamma$-comodules where $I_{j}$ is an injective comodule with no
$v_{n}$-torsion.  In particular, $L_{n}I_{j}=I_{j}$ by
Theorem~\ref{thm-localization}(a).  If we apply $L_{n}$ to this
sequence, we find that $L_{n}C_{j+1}=C_{j+1}$ and $L_{n}^{i}C_{j+1}=0$
for $i>0$, completing the induction step.  
\end{proof}


\providecommand{\bysame}{\leavevmode\hbox to3em{\hrulefill}\thinspace}
\providecommand{\MR}{\relax\ifhmode\unskip\space\fi MR }
\providecommand{\MRhref}[2]{%
  \href{http://www.ams.org/mathscinet-getitem?mr=#1}{#2}
}
\providecommand{\href}[2]{#2}

\end{document}